\documentclass{icmart}

\contact[d.kuhn@bham.ac.uk]{
Daniela K\"uhn,
School of Mathematics,
University of Birmingham,
Edgbaston,
Birmingham,
B15 2TT,
UK}

\contact[d.osthus@bham.ac.uk]{
Deryk Osthus, 
School of Mathematics,
University of Birmingham,
Edgbaston,
Birmingham,
B15 2TT,
UK}

\usepackage{amssymb}
\usepackage{amsmath}
\usepackage{graphicx}
\usepackage{fancyhdr}
\usepackage{pstricks,tikz}

\newcommand{\eps}{\varepsilon} 

\newcommand{\cH}{\mathcal{H}}

\newtheorem{firstthm}{Proposition}[section]
\newtheorem{theorem}[firstthm]{Theorem}
\newtheorem{question}[firstthm]{Question}

\newtheorem{cor}[firstthm]{Corollary}
\newtheorem{problem}[firstthm]{Problem}
\newtheorem{conj}[firstthm]{Conjecture}

\newdimen\margin   
\def\textno#1&#2\par{%
   \margin=\hsize
   \advance\margin by -4\parindent
          \setbox1=\hbox{\sl#1}%
   \ifdim\wd1 < \margin
      $$\box1\eqno#2$$%
   \else
      \bigbreak
      \hbox to \hsize{\indent$\vcenter{\advance\hsize by -3\parindent
      \it\noindent#1}\hfil#2$}%
      \bigbreak
   \fi}

\def\COMMENT#1{}
\def\TASK#1{}

\def\eps{{\varepsilon}}
\newcommand{\ex}{\mathbb{E}}

\newcommand{\cC}{\mathcal{C}}
\newcommand{\cP}{\mathcal{P}}

\title[Hamilton cycles in graphs and hypergraphs: an extremal perspective]{Hamilton cycles in graphs and hypergraphs: an extremal perspective}

\author[Daniela K\"uhn and Deryk Osthus]

\begin{document}

\begin{abstract}
As one of the most fundamental and well-known NP-complete problems, the Hamilton cycle problem has been the subject of
intensive research.
Recent developments in the area have highlighted the crucial role played by the notions of expansion and quasi-randomness.
These concepts and other recent techniques have led to the solution of several long-standing problems in the area.
New aspects have also emerged, such as resilience, robustness and the study of Hamilton cycles in hypergraphs.
We survey these developments and highlight open problems, with an emphasis on extremal and probabilistic approaches. 
\end{abstract}



\maketitle


\section{Introduction} \label{intro}

A Hamilton cycle in a graph $G$ is a cycle that contains all the vertices of $G$.
The decision problem of whether a graph contains a Hamilton cycle is among Karp's original list of 
NP-complete problems~\cite{Karp}. 
Together with the satisfiability problem SAT and graph colouring, it is probably one of the 
most well-studied NP-complete problems.
The techniques and insights developed for these fundamental problems have also found applications to many more related and
seemingly more complex questions.

The main approach to the Hamilton cycle problem has been to prove natural sufficient conditions
which are best possible in some sense.
This is exemplified by Dirac's classical theorem~\cite{Dirac}:
\emph{every graph $G$ on $n \ge 3$ vertices whose minimum degree is at least $n/2$ contains a Hamilton cycle.}
More generally, one can ask the following `extremal' question:
what value of some easily computable parameter (such as the minimum degree) ensures the existence  of a Hamilton cycle?
The field has an enormous literature, so we concentrate on recent developments:
several long-standing conjectures have recently been solved and new techniques have emerged.
In particular, recent trends include the increasing role of probabilistic techniques and viewpoints
as well as approaches based on quasi-randomness.

Correspondingly, in this survey we will focus on the following topics:
regular graphs and expansion; optimal packings of Hamilton cycles and Hamilton decompositions; 
random graphs; uniform hypergraphs; counting Hamilton cycles.
Notable omissions include the following topics: Hamilton cycles with additional properties (e.g.~$k$-ordered Hamilton cycles); pancyclicity;
generalized degree conditions (e.g.~Ore- and Fan-type conditions); structural constraints (e.g.~claw-free and planar graphs)
as well as digraphs. Many results in these areas are covered e.g.~in the surveys by Gould~\cite{gouldsurvey, gouldsurveyIII} and Bondy~\cite{bondy}.
Digraphs are discussed in~\cite{KOsurvey}, though some very recent results on digraphs are also included here.

\section{Regular graphs and expansion}\label{sec:regexp}

\subsection{Dense regular graphs}

The union of two cliques as well as the complete almost balanced bipartite graph show that the minimum degree bound in
Dirac's theorem is best possible. The former graph is disconnected and the latter is not regular.
This led Bollob\'as~\cite{egt} as well as H\"aggkvist (see~\cite{jackson_for_Haggvist}) to (independently) make the following conjecture:
\emph{Every $t$-connected $d$-regular graph $G$ on $n$ vertices with $d \geq n/(t+1)$ is Hamiltonian.}%
\COMMENT{Bollob\'as's conjecture was stronger, with $D \geq n/(t+1)-1$.}
The case $t=2$ was settled in the affirmative by Jackson~\cite{jackson_for_Haggvist}.
\begin{theorem}[\cite{jackson_for_Haggvist}]\label{exact}
Every $2$-connected $d$-regular graph on $n$ vertices with $d \geq n/3$ is Hamiltonian.
\end{theorem}
However, Jung~\cite{jung} and independently Jackson, Li and Zhu~\cite{jlz}
gave a counterexample to the conjecture for $t \geq 4$.
Until recently, the only remaining case $t=3$ was wide open.
K\"uhn, Lo, Osthus and Staden~\cite{KLOS1, KLOS2} proved this case
for all large $n$.
\begin{theorem}[\cite{KLOS1, KLOS2}]\label{exact2}
There exists an integer $n_0$ such that every $3$-connected $d$-regular graph on $n \geq n_0$ vertices with $d \geq n/4$ is Hamiltonian.
\end{theorem}
The theorem is best possible in the sense that the bound on $d$ cannot be reduced and $3$-connectivity cannot be replaced by $2$-connectivity.
The key to the proof is a structural partition result for dense regular graphs which was proved recently by the same authors~\cite{KLOS1}:
the latter result gives a partition of an arbitrary dense regular graph into a small number of `robust components', with very few edges 
between these components. Each robust component is either a `robust expander' or a `bipartite robust expander'.
Here a graph $G$ is a robust expander if for every set $S\subseteq V(G)$ of `reasonable size', its neighbourhood $N(S)$
is significantly larger than $S$, even after
some vertices and edges of $G$ are deleted (the precise definition is given in Section~\ref{sec:expand}). \cite{KLOS1} also contains further applications of this partition result.
Similar ideas might also be useful to prove analogues of Theorem~\ref{exact} (say) for directed and oriented graphs
(see \cite{KOsurvey} for such conjectured analogues).

Christofides, Hladk\'y and M\'ath\'e~\cite{CHM} used an approach related to that in the proof of Theorem~\ref{exact2}
to prove the famous `Lov\'asz conjecture' in the case of dense graphs.
\begin{conj} \label{lovasz}
Every connected vertex-transitive graph has a Hamilton path.
\end{conj}
In contrast to common belief, Lov\'asz~\cite{Lovasz} in 1969 actually asked for the construction of a connected vertex-transitive
graph containing no Hamilton path. Traditionally however, the Lov\'asz conjecture is always
stated in the positive.
A related folklore conjecture is the following:
\begin{conj} \label{cayley}
Every connected Cayley graph on at least three vertices contains a Hamilton cycle.
\end{conj}
Here a Cayley graph is defined as follows:
Let $H$ be a finite group and let $S \subseteq H$ be a subset with $S=S^{-1}$ such that $S$ does
not contain the identity.%
    \COMMENT{these conditions ensure that we may treat $C(H; S)$ as an undirected graph which does not
contain any loops}
The corresponding
Cayley graph $G(H; S)$ has vertex set equal to $H$. Two vertices $g, h  \in H$ are joined by a
edge if and only if there exists $s \in S$ such that $g = sh$. (So every Cayley graph is vertex-transitive.)

Maru\v{s}ic~\cite{Marusic} proved Conjecture~\ref{cayley} in the case when $H$ is abelian.
Alspach~\cite{alspach} conjectured that in this case one even obtains a decomposition
of the set of edges of $G(H; S)$ into edge-disjoint Hamilton cycles and at most one perfect matching.
For a survey of results on these conjectures, see for example~\cite{Kutnar}.

The following result of Christofides, Hladk\'y and M\'ath\'e~\cite{CHM} confirms the `dense' case of both Conjecture~\ref{lovasz} and~\ref{cayley}.
\begin{theorem}[\cite{CHM}]
For every $\eps > 0$ there exists an integer $n_0$ such that every connected vertex-transitive graph on
$n \ge n_0$ vertices of degree at least $\eps n$ contains a Hamilton cycle.
\end{theorem}
To prove this result, Christofides, Hladk\'y and M\'ath\'e define the notion of `iron-connectedness' which is related to that of robust expansion and consider a partition of the given vertex-transitive graph into `iron-connected' components.
It would be interesting to find out whether such a partition-based approach can also be extended to sparser graphs.

\subsection{Sparse graphs: Toughness and expansion}
The extremal examples for Theorem~\ref{exact2}
indicate that an obstacle to the existence of a Hamilton cycle is the fact that the graph is `easy to separate'
into several pieces. 
The examples also show that connectivity is not the appropriate notion to use in this context.
So a fruitful direction of research has been to study notions which are stronger than connectivity.


One of the most famous conjectures in this
direction is the toughness conjecture of Chv\'atal~\cite{Chvatal}. 
It states that if a graph is `hard to separate' into many pieces, then it contains a Hamilton cycle.
\begin{conj}[\cite{Chvatal}]\label{toughness}
There is a constant $t$
so that every $t$-tough graph has a 
Hamilton cycle. 
\end{conj}
Here a graph is \emph{$t$-tough} if, for every nonempty set $S\subseteq V(G)$, the graph $G-S$
has at most $|S|/t$ components. Trivially, every graph with a Hamilton cycle is $1$-tough.
Little progress has been made on this conjecture -- we only know that if the conjecture holds, then we must have $t \ge 9/4$~\cite{BBV}.

So instead of considering toughness, it has been more rewarding to consider the related (and in some sense stronger) notions of expansion
and quasi-randomness. 
By expansion, we usually mean the following: 
every small set $S$ of vertices has a neighbourhood $N(S)$ which is large compared to $|S|$
(more formally, $N(S)$ denotes the set of all those vertices which are adjacent to at least one vertex in $S$).
It is well known that expansion is closely linked to eigenvalues of the adjacency matrix:
a large eigenvalue gap is equivalent to good expansion properties (in which case we often call such a graph quasi-random).
In particular, there is a conjecture of Krivelevich and Sudakov~\cite{KS} on Hamilton cycles in regular graphs which
involves the `eigenvalue gap'. 
The conjecture itself would follow from the toughness conjecture.
\begin{conj}[\cite{KS}] 
There is a constant $C$ such that whenever
$G$ is a $d$-regular graph and the second largest (in absolute value) eigenvalue of the adjacency matrix of $G$ is a most
$d/C$, then $G$ has a Hamilton cycle.
\end{conj}

The best result towards this was proved by Krivelevich and Sudakov~\cite{KS}.
\begin{theorem}[\cite{KS}]\label{KSexpand}
There exists an integer $n_0$ such that the following holds for all $n\ge n_0$.
Suppose that $G$ is a $d$-regular graph on $n$ vertices and that the second largest (in absolute value)
eigenvalue $\lambda$ of the adjacency matrix of $G$ 
satisfies 
$\lambda  \le \frac{(\log \log n)^2}{1000 \log n (\log \log \log n)} d.
$
Then $G$ has a Hamilton cycle.
\end{theorem}
It is known that $\lambda=\Omega(d^{1/2})$ for $d\le n/2$. So the above result applies for example to quasi-random graphs
with $\lambda=\Theta(d^{1/2})$ whose density is polylogarithmic in~$n$, i.e.~for quasi-random graphs which are quite sparse.

The proof of Theorem~\ref{KSexpand} makes crucial use of the fact that the eigenvalue condition implies the following:
small sets of vertices expand 
and there are edges between any two large sets of vertices. Hefetz, Krivelevich and Szab\'o~\cite{HKS} proved the following general result
which goes beyond the class of regular graphs and makes explicit use of these conditions.
\begin{theorem}[\cite{HKS}]
There exists an integer $n_0$ such that the following holds for all integers $n, d$ with $n\ge n_0$ and
$12 \le d \le e^{(\log n)^{1/2}}$.
Let $m:=\frac{ n (\log \log n) \log d}{d \log n \log \log \log n}$. Suppose that
$G$ is a graph on $n$ vertices such that
$|N(S)| \ge  d|S|$ for every $S \subseteq V(G)$ with $|S| \le  m$. Moreover, suppose that
there is an edge in $G$ between any two disjoint subsets $A,B \in V(G)$ with $|A|,|B| \ge m/4130$.
Then $G$ has a Hamilton cycle.
\end{theorem}
The original motivation for this result was a problem on maker-breaker games, but the result also has several
other applications, see~\cite{HKS}.%
\COMMENT{this has applications to Hamiltonicity of Random Cayley graphs. But best result is due to Christofides and Markstrom}


\section{Packings of Hamilton cycles and decompositions}

\subsection{Optimal packings of Hamilton cycles in dense graphs}\label{sec:HCpack}

Nash-Williams~\cite{Diracext} proved a striking extension of Dirac's theorem: 
every graph on $n\ge 3$ vertices with minimum degree at least $n/2$ contains not just one but at least $5n/224$ edge-disjoint Hamilton cycles.
He conjectured~\cite{initconj,Diracext} that there should even be $n/4$ of these.
This was disproved by Babai (see~\cite{initconj}), who gave a construction showing that one cannot hope for more than
(roughly) $n/8$ edge-disjoint Hamilton cycles (see below for details).
Nash-Williams subsequently
raised the question of finding the best possible bound, which is answered in Corollary~\ref{NWmindegcor} below.

Recently Csaba, K\"uhn, Lapinskas, Lo, Osthus and Treglown~\cite{paper1,paper2,paper3,KLOmindeg} were able to answer a more general form of this question: 
\emph{what is the maximum number of edge-disjoint Hamilton cycles one can guarantee in a graph $G$ of minimum degree~$\delta$?}

A natural upper bound is obtained by considering the largest degree $\textnormal{reg}_{\textnormal{even}}(G)$
of an even-regular spanning subgraph of $G$. Let
\[
\textnormal{reg}_{\textnormal{even}}(n,\delta):=\min\{\textnormal{reg}_{\textnormal{even}}(G):|V(G)|=n,\ \delta(G)=\delta\}.
\]
Clearly, in general we cannot guarantee more than $\textnormal{reg}_{\textnormal{even}}(n,\delta)/2$
edge-disjoint Hamilton cycles in a graph of order $n$ and minimum
degree $\delta$. The next result of Csaba, K\"uhn, Lapinskas, Lo, Osthus and Treglown~\cite{paper1,paper2,paper3,KLOmindeg}
shows that this bound is best possible (if $\delta < n/2$, then $\textnormal{reg}_{\textnormal{even}}(n,\delta)=0$).%

\begin{theorem}[\cite{paper1,paper2,paper3,KLOmindeg}]\label{NWmindeg}
There exists an integer $n_0$ such that every graph $G$ on $n\ge n_0$ vertices contains at least
${\rm reg}_{\rm even}(n,\delta)/2$ edge-disjoint Hamilton cycles.
\end{theorem}

The main result in~\cite{KLOmindeg} proves Theorem~\ref{NWmindeg} unless $G$ is close to one of the two extremal graphs for
Dirac's theorem.
This allows us in~\cite{paper1,paper2,paper3} to restrict our attention to the latter situation (i.e.~when $G$ is close to
the complete balanced bipartite graph or close to the union of two disjoint copies of a clique).

An approximate version of Theorem~\ref{NWmindeg} for $\delta \ge n/2+\eps n$
was obtained earlier by Christofides, K\"uhn and Osthus~\cite{CKO}.
Hartke and Seacrest~\cite{HartkeHCs} gave a simpler argument with improved error bounds.

The parameter ${\rm reg}_{\rm even}(n,\delta)$ can be evaluated via Tutte's theorem. It turns out that for $n/2\le\delta<n$,
we have 
$$
\textnormal{reg}_{\textnormal{even}}(n,\delta) \sim
\frac{\delta+\sqrt{n(2\delta-n)}}{2},
$$
(see~\cite{CKO,Hartkefactors}).
In particular, if $\delta \geq n/2$ then ${\rm reg}_{\rm even}(n,\delta) \ge (n-2)/4$.
So Theorem~\ref{NWmindeg} implies the following explicit bound, which is best possible and
answers the above question of Nash-Williams~\cite{ initconj, Diracext}.
\begin{cor}\label{NWmindegcor}
There exists an integer $n_0$ such that every graph $G$ on $n\ge n_0$ vertices with
minimum degree $\delta(G)\ge n/2$ contains at least $(n-2)/8$ edge-disjoint Hamilton cycles.
\end{cor}
The following construction (which is based on a construction of Babai, see~\cite{initconj}) shows that the bound in 
Corollary~\ref{NWmindegcor} is best possible for $n=8k+2$, where $k \in \mathbb{N}$. 
Consider the graph $G$ consisting of one empty
vertex class $A$ of size $4k$, one vertex class $B$ of size $4k+2$ containing
a perfect matching and no other edges, and all possible edges between $A$ and $B$.
Thus $G$ has order $n=8k+2$ and minimum degree $4k+1=n/2$. 
Any Hamilton cycle in $G$ must contain at least two edges of
the perfect matching in~$B$, so $G$ contains at most $\lfloor |B|/4 \rfloor =k=(n-2)/8$ edge-disjoint
Hamilton cycles. 

A weaker version of Theorem~\ref{NWmindeg} for digraphs was proved by K\"uhn and Osthus in~\cite{KellyII}.
Ferber, Krivelevich and Sudakov~\cite{FKS} asked whether one can also obtain such a result for oriented graphs.

Recall that Theorem~\ref{NWmindeg} is best possible for the class of graphs on $n$ vertices with minimum degree $\delta$.
The following conjecture  of 
K\"uhn, Lapinskas and Osthus~\cite{KLOmindeg} would strengthen this in the sense that it would be best possible for every single graph $G$.
\begin{conj}[\cite{KLOmindeg}] \label{con:betterconj}
Suppose that $G$ is a graph on $n$ vertices
with minimum degree $\delta(G) \ge n/2$. Then $G$ contains
$\textnormal{reg}_{\textnormal{even}}(G)/2$ edge-disjoint
Hamilton cycles.
\end{conj}
For $\delta\ge(2-\sqrt{2}+\eps)n$, this conjecture was proved by K\"uhn and Osthus~\cite{KellyII}.
Recently, Ferber, Krivelevich and Sudakov~\cite{FKS} were able to obtain an approximate version of Conjecture~\ref{con:betterconj},
i.e.~a set of $(1-\eps)\textnormal{reg}_{\textnormal{even}}(G)/2$ edge-disjoint Hamilton cycles under the assumption that
$\delta(G) \ge (1+\eps)n/2$.

Also, it seems that the following `dual' version of the problem has not been investigated yet.
\begin{question}\label{q:cover}
Given a graph $G$ on $n$ vertices with $\delta(G)>n/2$, how many Hamilton cycles are needed in order to cover all the edges of $G$?
\end{question}
A trivial lower bound would be given by $\lceil \Delta(G)/2\rceil$.
However, this cannot always be achieved. Indeed, consider for example the graph $G$ obtained from a complete graph
on an odd number $n$ of vertices by deleting an edge~$xy$. Let $\mathcal{C}$ be a collection of Hamilton cycles covering all edges of~$G$.
Since both $x$ and $y$ have odd degree, at least one edge at each of $x$ and
$y$ has to lie in at least two Hamilton cycles from~$\mathcal{C}$. Thus $|\mathcal{C}|>(n-1)/2=\Delta(G)/2$.

Moreover, it is easy to see that the condition that $\delta> n/2$ in Question~\ref{q:cover}
is needed to ensure that every edge lies in a Hamilton cycle
(consider the balanced complete bipartite graph with a single edge in one of the classes).
More is known about the probabilistic version of Question~\ref{q:cover} (see Section~\ref{sec:random}).

Question~\ref{q:cover} can be viewed as  a restricted version of  
the following conjecture of Bondy~\cite{bondy90}, where arbitrary cycle lengths are permitted:%
\COMMENT{there is a related conjecture by Bermond et al in the handbook which says that cycles with a total length of $7e/5$ are enough}
\begin{conj}[\cite{bondy90}]
The edges of every $2$-edge-connected graph on $n$ vertices
can be covered by at most $2(n-1)/3$ cycles.%
\end{conj}


\subsection{The Hamilton decomposition and $1$-factorization conjectures}\label{sec:Hdec_1fac}

Theorem~\ref{NWmindeg} shows that for dense graphs the bottleneck for finding many edge-disjoint Hamilton cycles is the 
densest even-regular spanning subgraph. 
This makes it natural to consider the class of dense regular graphs.
In fact, Nash-Williams~\cite{initconj} suggested that these should even have a Hamilton decomposition.

Here a \emph{Hamilton decomposition} of a graph~$G$ consists of a set of edge-disjoint Hamilton cycles covering all edges of~$G$.
A natural extension of this to regular graphs $G$ of odd degree is to ask for a decomposition into 
Hamilton cycles and one perfect matching (i.e.~one perfect matching $M$ in $G$ together with a 
Hamilton decomposition of $G-M$). 
The most basic result in this direction is Walecki's theorem (see~\cite{lucas}), which dates back to the 19th century:
\begin{theorem}[see~\cite{lucas}]\label{Walecki}
If $n$ is odd, then the complete graph $K_n$ on $n$ vertices has a Hamilton decomposition.
If $n$ is even, then $K_n$ has a decomposition into Hamilton cycles together with a perfect matching.
\end{theorem}

The following result of Csaba, K\"uhn, Lo, Osthus and Treglown~\cite{paper1,paper2,paper3,paper4}
generalizes Walecki's theorem to arbitrary regular graphs which are sufficiently dense:
it determines the degree threshold for a regular graph to have a Hamilton decomposition.
In particular, it solves the above `Hamilton decomposition conjecture' of Nash-Williams~\cite{initconj} for all large graphs.

\begin{theorem}[\cite{paper1,paper2,paper3,paper4}] \label{HCDthm} 
There exists an integer $n_0$ such that the following holds.
Let $ n,d \in \mathbb N$ be such that $n \geq n_0$ and
$d \ge   \lfloor n/2 \rfloor $.
Then every $d$-regular graph $G$ on $n$ vertices has a decomposition into Hamilton cycles and 
at most one perfect matching.
\end{theorem}
The bound on the degree in Theorem~\ref{HCDthm} is best possible. 
Indeed, it is easy to see that a smaller degree bound would not even ensure connectivity.
Previous results include the following:
Nash-Williams~\cite{NWreg} showed that the degree bound in Theorem~\ref{HCDthm} guarantees a single Hamilton cycle.
Jackson~\cite{Jackson79} showed that one can guarantee close to $d/2-n/6$ edge-disjoint Hamilton cycles.
Christofides, K\"uhn and Osthus~\cite{CKO} obtained an approximate decomposition under the assumption that $d \ge n/2 +\eps n$.
Under the same assumption, K\"uhn and Osthus~\cite{KellyII} obtained an exact decomposition
(as a consequence of Theorem~\ref{decomp} below). 
Note that Conjecture~\ref{con:betterconj} would `almost' imply Theorem~\ref{HCDthm}.

Theorem~\ref{HCDthm} is related to the so-called `1-factorization conjecture'.
Recall that Vizing's theorem states that for any graph~$G$ of maximum degree $\Delta(G)$, the edge-chromatic number
$\chi'(G)$ of $G$ is either $\Delta(G)$ or $\Delta(G)+1$. 
For regular graphs $G$, $\chi'(G)=\Delta(G)$ is equivalent to the existence of a \emph{$1$-factorization},
i.e.~of a set of edge-disjoint perfect matchings covering all edges of~$G$.
The long-standing $1$-factorization conjecture guarantees a $1$-factorization in every regular graph of sufficiently high degree.
It was first stated explicitly by Chetwynd and Hilton~\cite{1factorization,CH} (who also proved partial results).
However, they state that
according to Dirac, it was already discussed in the 1950s. The following result of
Csaba, K\"uhn, Lo, Osthus and Treglown~\cite{paper1,paper2,paper3,paper4} confirms this conjecture for
sufficiently large graphs.

\begin{theorem}[\cite{paper1,paper2,paper3,paper4}]\label{1factthm}
There exists an $n_0$ such that the following holds.
Let $ n,d \in \mathbb N$ be such that $n\geq n_0$ is even and $d\geq 2\lceil n/4\rceil -1$. 
Then every $d$-regular graph $G$ on $n$ vertices has a $1$-factorization.
    Equivalently, $\chi'(G)=d$.
\end{theorem}
The bound on the minimum degree in Theorem~\ref{1factthm} is best possible. Indeed, a smaller bound on $d$ would
not even ensure a single perfect matching. To see this, suppose for example that $n = 2 \mod 4$
and consider the graph which is the disjoint union of two cliques of order $n/2$ (which is odd). 

Note that Theorem~\ref{HCDthm} does not quite imply Theorem~\ref{1factthm},
as the degree threshold in the former result is slightly higher.
The $1$-factorization conjecture is a special case of the `overfull subgraph' conjecture.
This would give an even wider class of graphs whose edge-chromatic number equals the maximum degree (see e.g.~the monograph~\cite{stiebitz}).

The best previous result towards the $1$-factorization conjecture is due to Perkovic and Reed~\cite{reed},
who proved an approximate version,
i.e.~they assumed that $d \ge n/2+\eps n$.
This was generalized by Vaughan~\cite{vaughan} to multigraphs of bounded multiplicity.

The following `perfect 1-factorization conjecture' was posed by Kotzig~\cite{Kotzig} more than fifty years ago at the first international conference devoted to Graph Theory.
It combines $1$-factorizations and Hamilton decompositions.
First note that it is easy to see that the complete graph $K_{2n}$ has a $1$-factorization.
The `perfect 1-factorization conjecture' would provide a far-reaching generalization of this fact.
\begin{conj}[\cite{Kotzig}]
$K_{2n}$ has a perfect $1$-factorization, i.e.~a $1$-factorization in which
any two $1$-factors induce a Hamilton cycle. 
\end{conj}
The conjecture is known to hold if $n$ or $2n-1$ is a prime, and for several special values of $n$, 
but beyond that very little is known.
To approach the conjecture it would be interesting to find $1$-factorizations so that the number of pairs of $1$-factors
which induce Hamilton cycles is as large as possible (see e.g.~\cite{wagner}).

Walecki's theorem can also be generalized in another direction: Alspach conjectured that one can decompose the complete graph
$K_n$ into cycles of arbitrary length. This was recently confirmed by Bryant, Horsley and Pettersson~\cite{BHP}.
\begin{theorem}
$K_n$ has a  decomposition into $t$ cycles of specified lengths $m_1,\dots,m_t$ if and only if $n$ is odd, $3 \le m_i \le n$ for 
$i\le t$, and $m_1+ \dots+m_t=\binom{n}{2}$. 
\end{theorem}
Perhaps it might be possible to prove a probabilistic analogue of this or extend the result to non-complete graphs.

As the final open problem in the area, we turn to a beautiful conjecture of 
Bermond (see~\cite{BermondL(G)}) that the existence of a Hamilton decomposition in a graph is inherited by its line graph
(note that an Euler circuit in a graph corresponds to a Hamilton cycle in the line graph).

\begin{conj}[see~\cite{BermondL(G)}]\label{Bermondlinegr}
If $G$ has a Hamilton decomposition, then the line graph $L(G)$ of $G$ has a Hamilton decomposition as well.
\end{conj}
Muthusamy and Paulraja~\cite{MP} proved this conjecture in the case when the number of Hamilton cycles in a Hamilton 
decomposition of $G$ is even (i.e.~when $G$ is $d$-regular where $4|d$). They also came quite close to proving it in the remaining 
case: they showed that if the number of Hamilton cycles in a Hamilton 
decomposition of $G$ is odd, then $L(G)$ can be decomposed into Hamilton cycles and one 2-factor.

\subsection{Kelly's conjecture}\label{sec:kelly}
Kelly's conjecture (see e.g.~\cite{moon}) dates back to 1968 and states that every regular tournament has a Hamilton decomposition. 
So one could view this as an oriented version of Walecki's theorem.
K\"uhn and Osthus~\cite{Kelly} recently proved the following result, which shows that Kelly's conjecture
is even true if one replaces the class of regular tournaments by that of sufficiently dense regular oriented graphs.
(An \emph{oriented graph} $G$ is a directed graph without $2$-cycles. $G$ is \emph{$d$-regular}
if all the in- and outdegrees equal~$d$.)

\begin{theorem}[\cite{Kelly}]\label{orientcor}
For every $\eps>0$ there exists an integer $n_0$ such that every $d$-regular oriented graph $G$ on $n\ge n_0$
vertices with $d\ge 3n/8+\eps n$ has a Hamilton decomposition. 
\end{theorem}
In fact, K\"uhn and Osthus deduce this result from an even more general result, which involves
an expansion condition rather than a degree condition (see Theorem~\ref{decomp}).
It is not clear whether the bound `$3n/8$' is best possible. 
However, this bound is a natural barrier since the minimum in- and outdegree threshold
which guarantees a single Hamilton cycle in an (not necessarily regular) oriented graph is $(3n-4)/8$. 
As mentioned above, Theorem~\ref{orientcor} implies Kelly's conjecture for all large tournaments.
\begin{cor}
There exists an integer $n_0$ such that
every regular tournament on $n\ge n_0$ vertices has a Hamilton decomposition. 
\end{cor}
K\"uhn and Osthus~\cite{KellyII} also used Theorem~\ref{orientcor} to prove a conjecture of
Erd\H{o}s on optimal packings of Hamilton cycles in random tournaments, which can be viewed as a probabilistic version
of Kelly's conjecture:
\begin{theorem}[\cite{KellyII}]\label{erdos}
Let $T$ be a tournament on $n$ vertices which is chosen uniformly at random. Then a.a.s.~$T$ contains
$\min \{\delta^+(T),\delta^-(T) \}$ edge-disjoint Hamilton cycles.
\end{theorem}
(Here we write a.a.s. for `asymptotically almost surely', see Section~\ref{sec:random} for the definition.)
The bound is clearly best possible. A similar phenomenon has been shown to occur in the random graph $G_{n,p}$ (see Theorem~\ref{HamDecomp}).

Jackson~\cite{jackson} posed the following bipartite version of Kelly's conjecture.
Here a \emph{bipartite tournament} is an orientation of a complete bipartite graph.
\begin{conj}[\cite{jackson}] \label{kellybip}
Every regular bipartite tournament has a Hamilton decomposition.
\end{conj}

It is not even known whether there exists an approximate decomposition, 
i.e.~a set of Hamilton cycles covering almost all the edges of a regular bipartite tournament.
Another conjecture related to Kelly's conjecture was posed by Thomassen.
The idea is to force many edge-disjoint Hamilton cycles by high connectivity rather than regularity:
Thomassen~\cite{thomassenconj} conjectured that for every $k$ there is an integer $f(k)$ so
that every strongly $f(k)$-connected tournament contains $k$ edge-disjoint 
Hamilton cycles. K\"uhn, Lapinskas, Osthus and Patel~\cite{KLOP} proved this by showing that $f(k)=O(k^2 (\log k)^2)$
and conjectured that $f(k)=O(k^2)$.


\subsection{Robust expansion}\label{sec:expand}

As we already indicated in Section~\ref{sec:regexp}, there is an intimate connection between expansion and Hamiltonicity.
In what follows, we describe a relatively new `dense' notion of expansion, which has been extremely fruitful in studying
not just Hamilton cycles but also Hamilton decompositions and more general subgraph embeddings.

Roughly speaking, this notion of `robust expansion' is defined as follows:
for any set $S$ of vertices, its robust neighbourhood is the set of all those vertices which have many neighbours in~$S$.
A graph is a robust expander if for every set $S$ which is not too small and not too large, its robust neighbourhood is at least
a little larger than~$S$ itself. 

More precisely, let $0<\nu\le  \tau<1$. Given any graph~$G$ on $n$ vertices and
$S\subseteq V(G)$, the \emph{$\nu$-robust neighbourhood~$RN_{\nu,G}(S)$ of~$S$}
is the set of all those vertices~$x$ of~$G$ which have at least $\nu n$ neighbours
in~$S$. $G$ is called a \emph{robust $(\nu,\tau)$-expander}
if 
$$
|RN_{\nu,G}(S)|\ge |S|+\nu n \ \mbox{ for all } \
S\subseteq V(G) \ \mbox{ with } \ \tau n\le |S|\le  (1-\tau)n.
$$
This notion was introduced (for digraphs) by K\"uhn, Osthus and Treglown~\cite{KOTchvatal}, who showed that every 
robustly expanding digraph of linear minimum in- and outdegree contains a Hamilton cycle.
Examples of robust expanders include graphs on $n$ vertices with minimum degree at least $n/2+\eps n$
as well as quasi-random graphs.
K\"uhn and Osthus~\cite{Kelly,KellyII} showed that every sufficiently large regular robust expander of linear
degree has a Hamilton decomposition.

\begin{theorem}[\cite{Kelly,KellyII}]\label{decomp}
For every $ \alpha >0$ there exists $\tau>0$ such that for all $\nu>0$ there exists an integer $n_0=n_0 (\alpha,\nu,\tau)$ 
for which the following holds. Suppose that
$G$ is a $d$-regular graph on $n \ge n_0$ vertices, where $d\ge \alpha n$, and that
$G$ is a robust $(\nu,\tau)$-expander.
Then $G$ has a Hamilton decomposition.
\end{theorem}
In~\cite{Kelly} they actually proved a version of this for digraphs, which has several applications.
(The undirected version is derived in~\cite{KellyII}.)
For example, this digraph version implies the following result.
\begin{theorem}[\cite{Kelly}]\label{regdigraph}
For every $\eps>0$ there exists an integer $n_0$ such that every $d$-regular digraph $G$ on $n\ge n_0$ vertices with $d\ge (1/2+\eps)n$
has a Hamilton decomposition.
\end{theorem} 
Theorem~\ref{regdigraph} is a far-reaching generalization of a result of Tillson, who proved a directed version of Walecki's theorem.
Moreover, Theorem~\ref{regdigraph} (which is algorithmic)
has an application to finding good tours for the (asymmetric) Traveling Salesman Problem (see~\cite{Kelly}).

The main original motivation for these results was to prove Kelly's conjecture for large tournaments:
indeed the directed version of Theorem~\ref{decomp} easily implies Theorem~\ref{orientcor}.

Theorem~\ref{decomp} has numerous further applications apart from Theorems~\ref{regdigraph} and~\ref{orientcor}
(both immediate ones and ones for which it is used as a tool).
For example, it is easy to see that for dense graphs, robust expansion is a relaxation of the notion of quasi-randomness.
So in particular, Theorem~\ref{decomp} implies (for large $n$) a recent result of Alspach, Bryant and Dyer~\cite{Paley}
that every Paley graph has a Hamilton decomposition. Theorem~\ref{decomp} is also used in the proof of the Hamilton
decomposition conjecture and the $1$-factorization conjecture
(Theorems~\ref{HCDthm} and~\ref{1factthm}).

The proof of Theorem~\ref{decomp} uses an `approximate' version of the result, which was proved by Osthus and Staden~\cite{OS}
and states that the conditions of the theorem imply the existence of an `approximate decomposition', 
i.e.~the existence of a set of edge-disjoint Hamilton cycles covering almost all edges of~$G$.
(This generalizes an earlier result of K\"uhn, Osthus and Treglown~\cite{KOTkelly} on approximate
Hamilton decompositions of regular tournaments.)


\section{Random graphs}\label{sec:random}

Probabilistic versions of the above Hamiltonicity questions have also been studied intensively. 
As usual, $G_{n,p}$ will denote a binomial random graph on $n$ vertices where every edge is present with probability $p$
(independently from all other edges),
and we say that a property of a random graph on $n$ vertices holds \emph{a.a.s.} (asymptotically almost surely) if the probability 
that it holds tends to 1 as $n$ tends to infinity.

Improving on bounds by several authors, Bollob\'as~\cite{Bol83}; Koml\'os and Szemer\'edi~\cite{KSz83} as well as Korshunov~\cite{Korshunov}
determined the precise value of $p$ which ensures a Hamilton cycle:
if $pn\ge \log n+ \log \log n +\omega(n)$, where $\omega(n) \to \infty$ as $n\to \infty$, then 
a.a.s.~$G_{n,p}$ contains a Hamilton cycle.
On the other hand, if $pn\le \log n+ \log \log n -\omega(n)$, then a.a.s.~$G_{n,p}$ contains an isolated vertex.

One can even obtain a `hitting time' version of this result in the evolutionary process $G_{n,t}$.
For this, let $G_{n,0}$ be the empty graph on $n$ vertices. Consider a random ordering of the edges of $K_n$. Let $G_{n,t}$
be obtained from $G_{n,t-1}$ by adding the $t$th edge in the ordering.
Ajtai, K\'omlos and Szemer\'edi~\cite{Komlos} as well as Bollob\'as~\cite{Bol84}
showed that a.a.s.~the time $t$ at which $G_{n,t}$ attains minimum degree two is the same as the time at which it 
first contains a Hamilton cycle.

There are many generalizations and related results. Recently, much attention has focused on
optimal packings of edge-disjoint Hamilton cycles and on resilience and robustness,
which we will discuss below. However, many intriguing questions remain open.


\subsection{Optimal packings of Hamilton cycles}

Bollob\'as and Frieze~\cite{BF85} extended the above hitting time result to packing edge-disjoint Hamilton cycles in random graphs of 
bounded minimum degree.
In particular, this implies the following: suppose that $pn \le \log n+ O(\log \log n)$.
Then a.a.s.~$G_{n, p}$ has
$\lfloor \delta(G_{n,p})/2 \rfloor$ edge-disjoint Hamilton cycles.
Frieze and Krivelevich~\cite{FK08} made the striking conjecture that this extends to all $p$.
This has recently been confirmed in a sequence of papers by several teams of authors:

\begin{theorem}\label{HamDecomp} For any $p=p(n)$, a.a.s.~$G_{n, p}$ has
$\lfloor \delta(G_{n,p})/2 \rfloor$ edge-disjoint Hamilton cycles.
\end{theorem}
We now summarize the results leading to a proof of Theorem~\ref{HamDecomp}.
Here `exact' refers to a bound of $\lfloor \delta(G_{n,p})/2 \rfloor$, `approx.' refers to
a bound of $(1-\eps)\delta(G_{n,p})/2$, and $\varepsilon$ is a positive constant.
$$
\begin{array}{l|l|l}
\mbox{authors}   & \mbox{range of } p &   \\
\hline
\mbox{Ajtai, Koml\'os, Szemer\'edi~\cite{Komlos}; Bollob\'as~\cite{Bol84}} &  \delta(G_{n,p}) =2 & \mbox{exact}  \\
\mbox{Bollob\'as \& Frieze~\cite{BF85}}   &  \delta(G_{n,p}) \mbox{ bounded} & \mbox{exact}  \\
\mbox{Frieze \& Krivelevich~\cite{fk}} & p \mbox{  constant} & \mbox{approx.}  \\
\mbox{Frieze \& Krivelevich~\cite{FK08}} & p=\frac{(1+o(1))\log n}{n} & \mbox{exact}  \\
\mbox{Knox, K\"uhn \& Osthus~\cite{AHDoRG}} & p \gg \frac{\log n}{n} & \mbox{approx.}  \\
\mbox{Ben-Shimon, Krivelevich \& Sudakov~\cite{BKS}} & \frac{(1+o(1))\log n}{n}\le p\le \frac{1.02\log n}{n} & \mbox{exact}  \\
\mbox{Knox, K\"uhn \& Osthus~\cite{KKO}} & \frac{(\log n)^{50}}{n} \le p \le 1- n^{-1/5} & \mbox{exact}  \\
\mbox{Krivelevich \& Samotij~\cite{KrS}} & \frac{\log n}{n} \le p \le n^{-1 + \eps} & \mbox{exact}  \\
\mbox{K\"uhn \& Osthus~\cite{KellyII}} & p \ge 2/3 & \mbox{exact}  \\
\end{array}
$$
In particular, the results in~\cite{BF85,KKO,KrS,KellyII} (of which~\cite{KKO,KrS} cover the main range)
together imply Theorem~\ref{HamDecomp}.

Glebov, Krivelevich and Szab\'o~\cite{GKS} were the first to consider the `dual' version of this problem:
\emph{how many Hamilton cycles are needed to cover all the edges of $G_{n,p}$?}
Hefetz, K\"uhn, Lapinskas and Osthus~\cite{HKLO} solved this problem for all $p$ that are not too small or too large
(based on the main lemma of~\cite{KKO}).
\begin{theorem}[\cite{HKLO}]\label{thm:main-result}
Suppose that $\frac{(\log n)^{117}}{n}\le p\le1-n^{-1/8}$.
Then a.a.s.~the edges of $G_{n,p}$ can be covered by $\lceil \Delta(G_{n,p})/2 \rceil$
Hamilton cycles.
\end{theorem}
It would be interesting to know whether a `hitting time' version of Theorem~\ref{thm:main-result} holds. For this, 
given a property $\cP$, let $t(\cP)$ denote the \emph{hitting time} of $\cP$, i.e.~the smallest $t$ so that $G_{n,t}$ has $\cP$. 
\begin{question}[\cite{HKLO}]
\label{con:hittime}
Let $\cC$ denote the property that an optimal covering of a graph $G$ with Hamilton cycles has size $\lceil \Delta(G)/2 \rceil$. 
Let $\cH$ denote the property that a graph $G$ has a Hamilton cycle.
Is it true that a.a.s.~$t(\cC)=t(\cH)$?
\end{question}
Note that $\cC$ is not monotone. In fact, it is not even the case that for all $t>t(\cC)$, $G_{n,t}$ a.a.s.~has $\cC$. Taking $n\ge 5$
odd and $t=\binom{n}{2}-1$, $G_{n,t}$ is the complete graph with one edge removed -- which, as noted at the end of Section~\ref{sec:HCpack},
cannot be covered by $(n-1)/2$ Hamilton cycles. It would be interesting to determine (approximately)
the ranges of $t$ such that a.a.s. $G_{n,t}$ has $\cC$.

Another natural model of random graphs is of course that of random regular graphs.
In this case it seems plausible that we can actually ask for a Hamilton decomposition
(and thus obtain an analogue of Theorem~\ref{HCDthm} for sparse random
graphs). Indeed, for random regular graphs of bounded degree this was proved by Kim and Wormald~\cite{KW01}
and for (quasi-)random regular graphs of linear degree this was proved by K\"uhn and Osthus~\cite{KellyII} (as a consequence of
Theorem~\ref{decomp}). However, the intermediate range remains open:
\begin{conj}
Suppose that $d=d(n) \to \infty$ and $d=o(n)$. Then a.a.s.~a random $d$-regular graph on $n$ vertices has a decomposition into 
Hamilton cycles and at most one perfect matching.
\end{conj}
So far, not even an approximate version of this is known. 
One might be able to deduce this from the results in~\cite{KKO}.

An analogue of the hitting time result of Bollob\'as and Frieze~\cite{BF85} for random geometric graphs was proved by
M\"uller, Perez-Gimenez and Wormald~\cite{MPW}.
Here the model is that $n$ vertices are placed at random on the unit square and edges are sequentially 
added in increasing order of edge-length. For fixed $k \ge 1$, they prove that a.a.s.~the first edge in the process that creates
minimum degree at least $k$ also causes the graph to have $\lfloor k/2 \rfloor$ edge-disjoint Hamilton cycles.
The hitting time result for the case $k=1$ was proved slightly earlier by
Balogh, Bollob\'as, Krivelevich, M\"uller and Walters~\cite{BBK}.


\subsection{Resilience}

Often one would like to know not just whether some graph $G$ has a property $\cP$, but `how 
strongly' it has this property.  In other words, does $G$ still have property $\cP$ if we delete (or add) some edges?
Implicitly, variants of this question have been studied for many properties and many classes of graphs.
Sudakov and Vu~\cite{SV} recently initiated the systematic study of this question.
In particular, they introduced the notion of \emph{resilience} of a graph with respect to a  property $\cP$
(below, we assume that $\cP$ is monotone increasing, i.e.~that $\cP$ cannot be destroyed by adding edges):
A graph has \emph{local resilience} $t$ with respect to $\cP$ if it still has $\cP$ whenever one deletes a set of edges
such that at each vertex less than $t$ 
edges are deleted. A graph has \emph{global resilience} $t$ with respect to $\cP$ if it still has $\cP$ whenever one deletes less than $t$ 
edges. 
Which of these variants is the more natural one to study usually depends on the property $\cP$:
for `global' properties such as Hamiltonicity and connectivity the local resilience leads to more interesting results, whereas
for `local' properties such as triangle containment, it makes more sense to study the global resilience.
%
Resilience has been studied intensively for various random graph models (mainly  $G_{n,p}$), as it yields natural probabilistic versions
of `classical' theorems. 
Lee and Sudakov~\cite{LS} proved a resilience version of Dirac's theorem (which improved previous bounds by several authors):
\begin{theorem}[\cite{LS}] 
For any $\eps>0$ there is a constant $C$ so that the following holds.
If $p \ge C \log n/n$ then a.a.s. every subgraph of $G_{n,p}$ with minimum degree at least $(1+\eps)np/2$
contains a Hamilton cycle.
\end{theorem}
It is natural to consider more general structures than Hamilton cycles.
However, as observed by Huang, Lee and Sudakov~\cite{HLS}, there is a limit to what one can ask for in this context:
for every $\eps>0$ there exists $p$ with $0<p<1$ such that a.a.s.~$G_{n,p}$ contains a subgraph $H$ with minimum degree at least $(1-\eps)np$
and $\Omega(1/p^2)$ vertices that are not contained in a triangle of~$H$.

As an even more informative notion than local resilience, Lee and Sudakov~\cite{LS} recently suggested a generalization of local resilience
which allows a different number of edges to be deleted at different vertices. In other words, in this `degree sequence resilience'
the degree sequence of the deleted graph has to be dominated by the given constraints. In particular, they asked for a
resilience version of Chv\'atal's theorem on Hamilton cycles:
\begin{problem}[\cite{LS}]
Characterize all those sequences $(k_1,\dots,k_n)$ for which $G=G_{n,p}$ a.a.s.~has the following property:
Let $H \subseteq G$ be such that the degree sequence $(d_1,\dots,d_n)$ of $H$ satisfies $d_i\le k_i$ for all $i\le n$.
Then $G-H$ has a Hamilton cycle.
\end{problem}
Partial results on this problem were obtained by Ben-Shimon, Krivelevich and Sudakov~\cite{BKS}.

\subsection{Robust Hamiltonicity}

An approach which can be viewed as `dual' to resilience was taken by Krivelevich, Lee and Sudakov~\cite{KLeeS}. 
They proved the following extension of Dirac's theorem, which one can view as a `robust' version of the theorem.
\begin{theorem}[\cite{KLeeS}] \label{robust}
There exists a constant $C$ such that for
$p \ge C \log n /n$ and a graph
$G$ on $n$ vertices of minimum degree at least $n/2$, the random subgraph
$G_p$ obtained from $G$ by including each edge with probability $p$ is a.a.s. Hamiltonian.
\end{theorem}
This theorem gives the correct order of magnitude of the threshold function since if
$p$ is a little smaller than $\log n/n$,
then the graph $G_p$ a.a.s.~has isolated vertices. Also, since there are graphs with minimum
degree $n/2-1$ which are not even connected, the minimum degree condition cannot be improved.
Note that the result can be viewed as an extension of Dirac's theorem since the case
$p = 1$ is equivalent to Dirac's theorem.

One can ask similar questions for other (families of) graphs which are known to be Hamiltonian.
In particular,
a natural question that seems to have been unfairly neglected is that of the Hamiltonicity threshold in random hypercubes.
More precisely, given $n$ and $p$, the random subgraph $Q_{n,p}$ of the $n$-dimensional cube~$Q_n$
is defined as follows: each edge of~$Q_n$ is included independently in $Q_{n,p}$ with probability
$p$. Bollob\'as~\cite{bollobasmatch} proved that if $p>1/2$ is a constant, then a.a.s. $Q_{n,p}$
is connected and has a perfect matching (and actually proved a hitting time version of this result).  
It seems plausible that a.a.s. $Q_{n,p}$ even contains a Hamilton cycle.
There is no chance for this if $p \le 1/2$ as there is a significant probability that $Q_{n,p}$ has an isolated 
vertex in that case.%
   \COMMENT{REFERENCE??}
\begin{conj}
Suppose that $p>1/2$ is a constant. Then a.a.s.~$Q_{n,p}$ has a Hamilton cycle.
\end{conj}
As far as we are aware, the question is still open even if $p$ is any constant close to one.
Since $Q_n$ is Hamiltonian, the above conjecture can be viewed as a `robust' version of this simple fact.%
  \COMMENT{subgraphs of quasi-random graphs - kriv}

\subsection{The P\'osa-Seymour conjecture}

Surprisingly, a probabilistic analogue of the P\'osa-Seymour conjecture is still open.
This beautiful generalization of Dirac's theorem 
states that every graph $G$ on $n$ vertices with minimum degree at least $kn/(k+1)$ contains the 
$k$th power of a Hamilton cycle
(which is obtained from a Hamilton cycle $C$ by adding edges between any vertices at distance at most $k$ on $C$). 
The conjecture was proved for large graphs by 
Koml\'os, S\'ark\"ozy and Szemer\'edi~\cite{KSSz98}. For squares of Hamilton cycles (i.e.~for $k=2$)
the best current bound in this direction is due to
Ch{\^a}u, DeBiasio and Kierstead~\cite{CBK}, who proved that in this case the conjecture
holds for all graphs on at least $2\cdot 10^8$ vertices.

A straightforward first moment argument indicates that the threshold for the square of a Hamilton cycle
in $G_{n,p}$ should be close to $p=n^{-1/2}$.
Note that unlike the deterministic version of the problem, this threshold would be significantly larger than 
the threshold for a triangle-factor. The latter was determined to be $n^{-2/3}(\log n)^{1/3}$ in a breakthrough by
Johansson, Kahn and Vu~\cite{JKV08}.

\begin{conj}[\cite{KOrandposa}] \label{square}
If $p \gg n^{-1/2}$, then a.a.s.~$G_{n,p}$ 
contains the square of a Hamilton cycle.
\end{conj}

When $k \ge 3$, the threshold is $n^{-1/k}$. 
This follows from a far more general theorem on thresholds for spanning structures in $G_{n,p}$ which was obtained by Riordan~\cite{riordan}.
His proof is based on the second moment method.
In~\cite{KOrandposa} K\"uhn and Osthus proved an `approximate' version of the above conjecture:
for any $\eps>0$, if $p \ge n^{-1/2+\eps}$, then $G_{n,p}$ a.a.s. contains the square of a Hamilton cycle.
Their proof is `combinatorial' in the sense that it uses a version of the absorbing method for random graphs rather than
the second moment method.%
\COMMENT{could ask for El-Zahar theorem for $G_{n,p}$, some work on this has been done by Frieze Kahn and Lubotzky as well as Kim et al.}
A version of this for quasi-random graphs was proved by Allen, B\"ottcher, H\`an, Kohayakawa and Person~\cite{ABHKP}.
Their result also extends to $k$th powers of Hamilton cycles.

In the spirit of Theorem~\ref{robust}, one could also ask about a `robust' version of Conjecture~\ref{square}.


\section{Hamilton cycles in uniform hypergraphs}

Cycles in hypergraphs have been studied since the 1970s.
The first notion of a hypergraph cycle was introduced by Berge~\cite{Berge}.
Recently, the much more structured notion of `$\ell$-cycles' has become very popular
and has led to very interesting results.

\subsection{Dirac-type theorems}

To obtain analogues of Dirac's theorem for hypergraphs, we first need to generalize the notions of a cycle and of
minimum degree. There are several natural notions available.

A \emph{$k$-uniform hypergraph} $G$ consists of a set $V(G)$ of vertices and a 
set $E(G)$ of edges so that each edge of consists of $k$ vertices. 
Given an integer $\ell$ with $1\le \ell<k$, we say that a $k$-uniform hypergraph~$C$ is an $\ell$-\textit{cycle} if there 
exists a cyclic ordering of the vertices of~$C$ such that every edge of~$C$ 
consists of~$k$ consecutive vertices and such that every pair of {consecutive} edges 
(in the natural ordering of the edges) intersects in precisely $\ell$ vertices. So every $\ell$-cycle $C$ has $|V(C)|/(k-\ell)$ edges.
In particular, $k-\ell$ divides the number of vertices in~$C$.
If $\ell=k-1$, then $C$ is called a \emph{tight cycle}, and if $\ell=1$, then $C$ is called a \emph{loose cycle}.
$C$ is a \emph{Hamilton $\ell$-cycle} of a $k$-uniform hypergraph~$G$ if $V(C)=V(G)$ and $E(C)\subseteq E(G)$.

More generally, a \emph{Berge cycle} is an alternating sequence $v_1,e_1,v_2,\dots,v_n,e_n$
of distinct vertices $v_i$ and distinct edges $e_i$ so that each $e_i$ contains $v_i$ and $v_{i+1}$.
(Here $v_{n+1}:=v_1$, and the edges $e_i$ are also allowed to contain vertices outside $\{v_1,\dots,v_n\}$.)
Thus every $\ell$-cycle is also a Berge cycle.
A Berge cycle $v_1,e_1,v_2,\dots,v_n, e_n$ is a \emph{Hamilton Berge cycle} of a hypergraph $G$ if $V(G)=\{v_1,\dots,v_n\}$
and $e_i\in E(G)$ for each $i\le n$. So a Hamilton Berge cycle of $G$ has $|V(G)|$ edges. Moreover, every
tight Hamilton cycle of $G$ is also a Hamilton Berge cycle of $G$ (but this is not true for Hamilton $\ell$-cycles with $\ell\le k-2$
as they have $|V(G)|/(k-\ell)$ edges).

We now introduce several notions of minimum degree for a $k$-uniform hypergraph~$G$.
Given a set $S$ of  vertices of $G$, the \emph{degree} $d_G(S)$ of~$S$ is the number of all those edges of $G$ 
which contain $S$ as a subset. The \emph{minimum $t$-degree} 
$\delta_t(G)$ of $G$ is then the minimum value of $d_G(S)$ taken over all 
sets $S$ of $t$ vertices of $G$. When $t=1$ we refer to this as the \emph{minimum vertex degree} of~$G$,
and when $t=k-1$ we refer to this as the \emph{minimum codegree}.

A Dirac-type theorem for Berge cycles was proved by Bermond, Germa, Heydemann and Sotteau~\cite{Bermond}.
A Dirac-type theorem for tight Hamilton cycles was proved by R\"odl, Ruci\'nski and Szemer\'edi~\cite{RRS,RRS2}.
(This improved an earlier bound by Katona and Kierstead~\cite{KK}.)
Together with the fact that if $(k-\ell)|n$ then any tight cycle contains an $\ell$-cycle on the same 
vertex set (consisting of every $(k-\ell)$th edge), this yields the following result.

\begin{theorem}[\cite{RRS,RRS2}]\label{rrs}
For all~$k\geq 3$, $1 \leq \ell \leq k-1$ and any~$\eps > 0$ there exists an integer~$n_0$ 
so that if $n\ge n_0$ and $(k-\ell)|n$ then any $k$-uniform hypergraph~$G$ on~$n$ vertices 
with $\delta_{k-1}(G) \geq \left(\frac{1}{2}+\eps\right) n$ contains a Hamilton 
$\ell$-cycle.
\end{theorem}

If $(k-\ell)|k$ and $k|n$ then the above result is asymptotically best possible. Indeed, to see this, note that if the above divisibility conditions
hold, then every $\ell$-cycle $C$ contains a perfect matching (consisting of every $k/(k-\ell)$th edge of $C$). On the other hand,
it is easy to see that the following parity based construction shows that a minimum codegree of $n/2-k$ does not
ensure a perfect matching: Given a set $V$ of $n$ vertices, let $A\subseteq V$ be a set of vertices such that $|A|$ is odd
and $n/2-1\le |A|\le n/2+1$. Let $G$ be the $k$-uniform hypergraph whose edges consists of all those $k$-element subsets $S$ of $V$
for which $|S \cap A|$ is even.

For $k=3$,  R\"odl, Ruci\'nski and Szemer\'edi~\cite{RRS3} were able to prove an exact version of Theorem~\ref{rrs}
(the threshold in this case is $\lfloor n/2\rfloor$).
The following result of K\"uhn, Mycroft and Osthus~\cite{KMO} deals with all those cases in which Theorem~\ref{rrs} is not
asymptotically best possible.

\begin{theorem}[\cite{KMO}]\label{main}
For all $k\geq 3$, $1 \leq \ell \leq k-1$ with $(k-\ell)\nmid k$ and any 
$\eps > 0$ there exists an integer $n_0$ so that if $n\ge n_0$ and $(k-\ell)|n$ then any 
$k$-uniform hypergraph $G$ on $n$ vertices with $$\delta_{k-1}(G) \geq \left(\frac{1}{\lceil 
\frac{k}{k-\ell} \rceil(k-\ell)}+\eps\right) n$$ contains a Hamilton 
$\ell$-cycle.
\end{theorem}
Theorem~\ref{main} is asymptotically best possible. To see this, let $t:=n/(k-\ell)$ and $s:=\lceil k/(k-\ell) \rceil$.
Fix a set $A$ of $\lceil t/s\rceil-1$ vertices and
consider the $k$-uniform hypergraph $G$ on $n$ vertices whose hyperedges all have nonempty intersection with $A$. Then $\delta_{k-1}(G)=|A|$.
However, an $\ell$-cycle on $n$ vertices has $t$ edges and every vertex on such a cycle lies in at most $s$
edges. So $G$ does not contain an Hamilton $\ell$-cycle since $A$ would be a vertex cover for such a cycle and $|A|s<t$.

So the problem of which codegree forces a Hamilton $\ell$-cycle is asymptotically solved,
though exact versions covering all cases remain a challenging open problem.
For $k=3$ and $\ell=1$, Czygrinow and Molla~\cite{CM} were able to prove such an exact version.
The following table describes the history of the results leading to the current state of the art.

$$
\begin{array}{l|l|l|l}
\mbox{authors}   &  k  & \mbox{range of } \ell &   \\
\hline
\mbox{R\"odl, Ruci\'nski \& Szemer\'edi~\cite{RRS}} & k=3 & \ell=2 & \mbox{approx.}  \\
\mbox{K\"uhn \& Osthus~\cite{KO}} & k=3 & \ell=1 & \mbox{approx.}  \\
\mbox{R\"odl, Ruci\'nski \& Szemer\'edi~\cite{RRS2}} & k\ge 3 & \ell=k-1 & \mbox{approx.}  \\
\mbox{Keevash, K\"uhn, Mycroft \& Osthus~\cite{KKMO}} & k\ge 3 & \ell=1 & \mbox{approx.}  \\
\mbox{H\`an \& Schacht~\cite{HS}} & k\ge 3 & 1\le \ell< k/2 & \mbox{approx.}  \\
\mbox{K\"uhn, Mycroft \& Osthus~\cite{KMO}} & k\ge 3 & 1\le \ell< k-1, \ (k-\ell)\nmid k & \mbox{approx.}  \\
\mbox{R\"odl, Ruci\'nski \& Szemer\'edi~\cite{RRS3}} & k= 3 & \ell=2 & \mbox{exact}  \\
\mbox{Czygrinow and Molla~\cite{CM}} & k= 3 & \ell=1 & \mbox{exact}  \\
\end{array}
$$

Proving corresponding results for vertex degrees seems to be considerably harder.
The following natural conjecture, which was implicitly posed by R\"odl and Ruci\'nski~\cite{RR}, is wide open.
\begin{conj}[\cite{RR}] \label{hypercycle}
For all integers $k\ge 3$ and all $\eps>0$ there is an integer $n_0$ so that the following holds: if~$G$ 
is a $k$-uniform hypergraph on $n\ge n_0$ vertices with
$$\delta_1(G)\ge \left(1-\left(1-\frac{1}{k}\right)^{k-1}+\eps\right)\binom{n}{k-1},$$ then $G$ contains a tight Hamilton cycle.
\end{conj}

This would be asymptotically best possible. Indeed, if $k|n$ then any tight Hamilton cycle contains a perfect matching,
and a minimum vertex degree which is slightly smaller than in Conjecture~\ref{hypercycle} would not even guarantee a perfect matching.
To see the latter, fix a set $A$ of $n/k-1$ vertices and
consider the $k$-uniform hypergraph $G$ on $n$ vertices whose hyperedges all have nonempty intersection with $A$.
Then $\delta_1(G)\sim (1-(1-1/k)^{k-1})\binom{n}{k-1}$, but $G$ does not contain a perfect matching.

For general $k$, Conjecture~\ref{hypercycle} seems currently out of reach -- it is even a major open question to determine whether the above degree
bound ensures a perfect matching of~$G$. However, it would also be interesting to obtain non-trivial bounds (see e.g.~\cite{RR}).
For $k=3$ the best current bound towards Conjecture~\ref{hypercycle} was proved by R\"odl and Ruci\'nski~\cite{RR2}. They showed
that in this case the conjecture holds if $1-(1-1/3)^{2}=5/9$ is replaced by $(5-\sqrt{5})/3$.

For $k=3$, Han and Zhao~\cite{HZ} were able to determine
the minimum vertex degree which guarantees a loose Hamilton cycle exactly.

\begin{theorem}[\cite{HZ}]
There exists an integer $n_0$ so that the following holds. Suppose that $G$ is a $3$-uniform hypergraph on $n\ge n_0$ vertices
with $\delta_1(G) \ge \binom{n}{2} -\binom{3n/4}{2}+c$, where $n$ is even, $c=2$ if $4|n$ and $c=1$ otherwise.
Then $G$ contains a loose Hamilton cycle.
\end{theorem}
The bound on the minimum vertex degree is tight:
for $n$ of the form $4t+2$, fix a set $A$ of $t$ vertices and
consider the $k$-uniform hypergraph $G$ on $n$ vertices whose hyperedges all have nonempty intersection with $A$.
Bu{\ss}, Han and Schacht~\cite{BHSch} had earlier proved an asymptotic version of this result.


\subsection{Hamilton cycles in random hypergraphs}
Similarly as in the graph case, it is natural to study Hamiltonicity questions in a probabilistic setting.
Let $H_{n,p}^{(k)}$ denote the random $k$-uniform hypergraph on $n$ vertices where every edge is present with probability~$p$,
independently of all other edges. The following result of Dudek, Frieze, Loh and Speiss~\cite{DFLS} determines the threshold for the
existence of a loose Hamilton cycle in~$H_{n,p}^{(k)}$.
(In both Theorems~\ref{loose} and~\ref{tight} we only consider those $n$ which satisfy the trivial divisibility condition
for the existence of an $\ell$-cycle, i.e.~that $n$ is a multiple of $k-\ell$.)

\begin{theorem}[\cite{DFLS}]\label{loose}
Suppose that $k \ge 3$. If $p \gg (\log n)/n^{k-1}$, then a.a.s.~$H_{n,p}^{(k)}$ contains a loose Hamilton cycle.
\end{theorem}
The logarithmic factor appears due to the `local' obstruction that a.a.s.~$H_{n,p}^{(k)}$ contains isolated 
vertices below this threshold.

The proof of Theorem~\ref{loose} is `combinatorial' (in particular, it does not use the second moment method)
and builds on earlier results by Frieze~\cite{Frieze_loose} as well as Dudek and Frieze~\cite{DFloose},
which required additional divisibility assumptions.
The argument in~\cite{DFLS} also uses the celebrated result of Johansson, Kahn and Vu~\cite{JKV08} on the threshold for perfect matchings in hypergraphs.%
   \COMMENT{for $k=3$ the result is slightly stronger: there exists a $c$ so that $p \ge c(\log n)/n^2$ is enough}

Loose Hamilton cycles in random regular hypergraphs have been considered by Dudek, Frieze, Ruci\'nski and  \v{S}ileikis~\cite{DFRS}.
The next result due to Dudek and Frieze~\cite{DFtight} concerns precisely those values of $k$ and $\ell$ not covered by Theorem~\ref{loose}.
Thus together Theorems~\ref{loose} and~\ref{tight} determine the threshold for the existence of a Hamilton $\ell$-cycle in
random $k$-uniform hypergraphs for any given value of~$k$ and~$\ell$.
\begin{theorem}[\cite{DFRS}]\label{tight}
\
\begin{itemize}
\item[{\rm (i)}] For all integers $k > \ell \ge 2$ and fixed $\eps > 0$, if $p = (1 - \eps)e^{k- \ell}/n^{k-\ell}$, then 
a.a.s.~$H^{(k)}_{n,p}$ does not contain a Hamilton 
$\ell$-cycle.
\item[{\rm (ii)}] If $k > \ell \ge 2$ and $p \gg 1 /n^{k-\ell}$, 
then a.a.s.~$H^{(k)}_{n,p}$ contains a Hamilton $\ell$-cycle.%
\COMMENT{ for $\ell\ge 3$ one actually gets that $p \ge c/n^{k-\ell}$ for large enough $c$ is enough}
\item[{\rm (iii)}] For all fixed $\eps > 0$, if $k \ge 4$ and $p = (1 + \eps)e/n$, then a.a.s.~$H^{(k)}_{n,p}$ contains a tight Hamilton cycle. 
\end{itemize}
\end{theorem}
The proof of Theorem~\ref{tight} is based on the second moment method (which seems to fail for Theorem~\ref{loose}).
An algorithmic proof of (iii) with a weaker threshold of $p \ge n^{-1+\eps}$ was given by Allen, B\"ottcher, Kohayakawa  and Person~\cite{ABKP}.
Note that, for $k\ge 4$, (i) and~(iii) establish a sharp threshold for tight Hamilton cycles, i.e.~when $\ell=k-1$. 
It would be interesting to obtain a sharp threshold for other cases besides those in~(iii)
and a hitting time result for loose Hamilton cycles.

\subsection{Hamilton decompositions}

Hypergraph generalisations of Walecki's theorem (Theorem~\ref{Walecki}) have also been investigated.
This question was first studied for the notion of a Berge cycle.
Let $K_n^{(k)}$ denote the complete $k$-uniform hypergraph on $n$ vertices. 
Since every Hamilton Berge cycle of $K_n^{(k)}$ has $n$ edges, a necessary condition for the existence of 
a decomposition of $K_n^{(k)}$ into Hamilton Berge cycles is that $n$ divides $\binom{n}{k}$.
Bermond, Germa, Heydemann and Sotteau~\cite{Bermond} conjectured that this condition is also sufficient.
For $k=3$, this conjecture follows by combining the results of Bermond~\cite{bermond} and Verrall~\cite{verrall}. 
K\"uhn and Osthus~\cite{KOBerge} showed that as long as $n$ is not too small, the conjecture holds for $k\ge 4$ as well.
So altogether this yields the following result.
\begin{theorem}[\cite{bermond, verrall, KOBerge}]\label{simple}
Suppose that $3 \le k <n$, that $n$ divides $\binom{n}{k}$ and, in the case when $k\ge 4$, that $n\ge 30$.
Then $K_n^{(k)}$ has a decomposition into Hamilton Berge cycles.
\end{theorem}

The following conjecture of K\"uhn and Osthus~\cite{KOBerge} would be an analogue of Theorem~\ref{simple} for Hamilton $\ell$-cycles.

\begin{conj}[\cite{KOBerge}]\label{lcycles}
For all integers $1\le \ell <k$ there exists an integer $n_0$ such that the following
holds for all $n\ge n_0$. Suppose that $k-\ell$ divides $n$ and that $n/(k-\ell)$ divides $\binom{n}{k}$.
Then $K_n^{(k)}$ has a decomposition into Hamilton $\ell$-cycles.
\end{conj}
To see that the divisibility conditions are necessary, recall that every $\ell$-cycle on $n$ vertices contains exactly $n/(k-\ell)$ edges.

The `tight' case $\ell=k-1$ of Conjecture~\ref{lcycles} was already formulated and investigated by Bailey and Stevens~\cite{baileystevens}.
Actually, if $n$ and $k$ are coprime, the case $\ell=k-1$ already corresponds to a conjecture made 
independently by Baranyai~\cite{bconj} and Katona concerning `wreath decompositions'.
A $k$-partite version of the `tight' case of Conjecture~\ref{lcycles} was recently proved by Schroeder~\cite{schroeder}.

Conjecture~\ref{lcycles} is known to hold `approximately' (with some additional divisibility conditions on $n$),
i.e.~one can find a set of edge-disjoint Hamilton $\ell$-cycles which together
cover almost all the edges of $K_n^{(k)}$. 
This is a very special case of results in~\cite{BF,FK,FKL} which together guarantee approximate decompositions of quasi-random
uniform hypergraphs into Hamilton $\ell$-cycles for $1 \le \ell < k$ (again, the proofs need $n$ to satisfy additional divisibility constraints).

For example, 
Frieze, Krivelevich and Loh~\cite{FKL} proved an approximate decomposition result for tight Hamilton cycles in quasi-random $3$-uniform hypergraphs, which implies the following
result about random hypergraphs. 
\begin{theorem}[\cite{FKL}]
Suppose that $\eps, p,n$ satisfy $\eps^{45}np^{16} \ge (\log n)^{21}$. Then whenever $4|n$, a.a.s.~there is a collection of edge-disjoint tight Hamilton cycles of $H_{n,p}^{(3)}$
which cover all but at most an $\eps^{1/15}$-fraction of the  edges of $H_{n,p}^{(3)}$.
\end{theorem}
The proof proceeds via a reduction to an approximate decomposition result of quasi-random digraphs into Hamilton cycles.
This reduction is also the cause for the divisibility requirement. 
It would be nice to be able to eliminate this requirement.
It would also be interesting to know whether the threshold for the existence of an approximate decomposition into Hamilton $\ell$-cycles 
coincides with the threshold for a single Hamilton cycle.

\section{Counting Hamilton cycles}

In Section~\ref{sec:HCpack} the aim was to strengthen Dirac's theorem (and other results) by finding many edge-disjoint Hamilton cycles.
Similarly, it is natural to omit the condition of edge-disjointness and ask for the total number of Hamilton cycles in a graph.
For Dirac graphs (i.e.~for graphs on $n$ vertices with minimum degree at least $n/2$), this problem was essentially solved by
Cuckler and Kahn~\cite{CucklerKahn,CucklerKahn2}.
They gave a remarkably elegant formula which asymptotically determines the logarithm of the number of Hamilton cycles.

To state their result, we need the following definitions.
For a graph $G$ and edge weighting $x \colon E(G) \rightarrow \mathbb{R}^+$, set 
$h(x) := \sum_{e\in E(G)} x_e \log_2(1/x_e)$, where $x_e$ denotes the weight of the edge~$e$. 
This is related to the entropy function, except that $\sum_{e\in E(G)} x_e$ is not required to equal 1. 
We call an edge weighting $x$ a \emph{perfect fractional matching} if
$\sum_{e \ni v} x_e = 1$ for each vertex $v$ of~$G$.  
Finally, let $h(G)$ (the `entropy' of $G$) be the maximum of
$h(x)$ over all fractional matchings $x$. 
\begin{theorem}[\cite{CucklerKahn,CucklerKahn2}]
Suppose that $G$ is a graph on $n$ vertices with $\delta(G) \ge n/2$.
Then the number of Hamilton cycles in $G$ is 
\begin{equation} \label{count1}
2^{2h(G) - n \log_2 e - o(n) }.
\end{equation}
In particular, the number of Hamilton cycles in $G$ is at least  
\begin{equation} \label{count2}
(1-o(1))^n\frac{\delta(G)^n}{n^n} n! \ge \frac{n!}{(2 + o(1))^n}.
\end{equation}
\end{theorem}
(\ref{count2}) answers a question of S\'ark\"ozy, Selkow and Szemer\'edi~\cite{SSS}.
The proof of the lower bound in~(\ref{count1}) proceeds by considering a random walk which embeds the Hamilton cycles.
(\ref{count2}) is a consequence of~(\ref{count1}), but the derivation is nontrivial.
(It is easy to derive if $G$ is $d$-regular, as then setting $x_e:=1/d$ for each edge $e$ of $G$ maximises $h(x)$.)
As a general bound on the number of Hamilton cycles in Dirac graphs, (\ref{count2}) is best possible (up to lower order terms)
-- consider for example the complete balanced bipartite graph.
In fact, it is an easy consequence of Bregman's theorem on permanents that the first bound in~(\ref{count2}) is best possible for
\emph{any} regular graph.

$h(G)$ can be computed in polynomial time, so one can efficiently obtain a rough estimate for the number
of Hamilton cycles in a given Dirac graph.
The question of obtaining more precise estimates via randomized algorithms was considered earlier by Dyer, Frieze and Jerrum~\cite{DFJ}.
For graphs whose minimum degree is at least  $n/2+\eps n$, they
obtained a fully polynomial time randomized approximation scheme (FPRAS) for counting the number of Hamilton cycles.
(Roughly speaking, an FPRAS is a randomized polynomial time algorithm which gives an answer to a counting problem
to within a factor of $1+o(1)$ with probability $1-o(1)$.)
They asked whether this result can be extended to all Dirac graphs.
\begin{question}[\cite{DFJ}]
Let $\mathcal{G}$ denote the class of all Dirac graphs, i.e.~of all graphs $G$ with minimum degree at least $|V(G)|/2$. 
Is there an FPRAS for counting the number of Hamilton cycles for all graphs in $\mathcal{G}$?
\end{question}

Ferber, Krivelevich and Sudakov~\cite{FKS} proved an analogue of~(\ref{count2}) for oriented graphs whose degree is slightly above the 
Hamiltonicity threshold.

Counting Hamilton cycles also yields interesting results in the random graph setting.
Note that the expected number of Hamilton cycles in $G_{n,p}$ is $p^n (n-1)!/2$.
Glebov and Krivelevich~\cite{GK} showed that for any $p$ above the Hamiltonicity threshold,
a.a.s.~the number of Hamilton cycles in $G_{n,p}$ is not too far from this.
\begin{theorem}[\cite{GK}]\label{countHC}
Let $p \ge \frac{\log n+\log \log n +\omega(n)}{n}$, where $\omega(n)$ tends to infinity with $n$. 
Then a.a.s.~the number of Hamilton cycles in $G_{n,p}$ is $(1 - o(1))^np^n n!$.
\end{theorem}
For $p = \Omega(n^{-1/2})$, this was already proved by Janson~\cite{Janson1}, who in fact determined the asymptotic distribution of the 
number of Hamilton cycles in this range. Surprisingly, his results imply that a.a.s.~the number $X$ of Hamilton cycles in $G_{n,p}$
is concentrated below the expected value, i.e.~a.a.s.~$X/\ex(X) \to 0$ for $p \to 0$ 
(on the other hand, in the $G_{n,m}$ model, $X$ is concentrated at $\ex(X)$ in the range when $n^{3/2}\le m\le 0.99\binom{n}{2}$).
Glebov and Krivelevich~\cite{GK} also obtained a hitting time version of Theorem~\ref{countHC}.
\begin{theorem}[\cite{GK}]\label{GKcounting} 
In the random graph process $G_{n,t}$, at the very moment the minimum degree
becomes two, a.a.s.~the number of Hamilton cycles becomes $(1 - o(1))^n(\log n /e)^n$.
\end{theorem}
Note that at the hitting time $t$ for minimum degree two a.a.s.~the edge density $p$ of $G_{n,t}$ is close to
$\log n/n$, and so the expression in Theorem~\ref{GKcounting} could also be written as $(1 - o(1))^np^nn!$,
which coincides with Theorem~\ref{countHC}.

A related result of Janson~\cite{Janson2} determines the asymptotic distribution of the number of Hamilton cycles in random $d$-regular graphs for constant~$d\ge 3$.
Frieze~\cite{Frcount} proved a similar formula to that in Theorem~\ref{countHC} for dense quasi-random graphs, which was extended to sparse quasi-random graphs by
Krivelevich~\cite{Krivelevich}.

It turns out that the number of Hamilton cycles in a graph is often closely connected to the number of perfect matchings
(indeed the former is always at most the square of the latter).
So most of the above papers also contain related results about counting perfect matchings.


\begin{thebibliography}{100}

\bibitem{Komlos} M.~Ajtai, J.~Koml\'os and E.~Szemer\'edi, The first occurrence of Hamilton cycles in random graphs, \emph{Ann. Discrete Math.}~{\bf 27} (1985), 173--178.


\bibitem{ABHKP} P. Allen, J.~B\"ottcher, H.~H\`an, Y.~Kohayakawa and Y.~Person, Powers of Hamilton cycles in pseudorandom graphs, preprint.

\bibitem{ABKP} P. Allen, J.~B\"ottcher, Y.~Kohayakawa and Y.~Person, 
Tight Hamilton cycles in random hypergraphs, \emph{Random Structures \& Algorithms}, to appear.

\bibitem{alspach} B. Alspach, Research Problem 59, \emph{Discrete Math.} {\bf 50} (1984), 115. 

\bibitem{BermondL(G)} B.~Alspach, J.C.~Bermond and D.~Sotteau,  Decomposition into cycles. I. Hamilton decompositions,
in: in: Cycles and rays (Montreal, PQ, 1987),  
\emph{NATO Adv. Sci. Inst. Ser. C Math. Phys. Sci.} {\bf 301}, Kluwer Acad. Publ., Dordrecht (1990), 9--18.

\bibitem{Paley} B.~Alspach, D.~Bryant and D.~Dyer, 
Paley graphs have Hamilton decompositions, \emph{Discrete Math.} {\bf 312} (2012), 113--118.

\bibitem{baileystevens} R.~Bailey and B.~Stevens,
Hamiltonian decompositions of complete $k$-uniform hypergraphs, \emph{Discrete Math.} {\bf 310} (2010), 3088--3095.

\bibitem{BF} D.~Bal and A.~Frieze,
Packing tight Hamilton cycles in uniform hypergraphs,
\emph{SIAM J. Discrete Math.} {\bf 26} (2012), 435--451.

\bibitem{BBK}
J. Balogh, B. Bollob\'as, M. Krivelevich, T. M\"uller and M. Walters, 
Hamilton cycles in random geometric graphs, 
\emph{Annals of Applied Probability} {\bf 21} (2011), 1053--1072. 

\bibitem{bconj} Zs. Baranyai, 
The edge-coloring of complete hypergraphs I,
\emph{J. Combin. Theory  B} {\bf 26} (1979), 276--294. 

\bibitem{BBV} D.~Bauer, H.J.~Broersma, H.J.~Veldman, Not every 2-tough graph is Hamiltonian,
\emph{Discrete Applied Math.}~\textbf{99} (2000), 317--321.



\bibitem{BKS} S.~Ben-Shimon, M.~Krivelevich and B.~Sudakov,
On the resilience of Hamiltonicity and optimal packing of Hamilton cycles in random graphs,
\emph{SIAM J. Discrete Math.} {\bf 25} (2011), 1176--1193.

\bibitem{Berge} C.~Berge, \emph{Graphs and Hypergraphs}, North-Holland, Amsterdam, 1979.

\bibitem{bermond} J.C.~Bermond, Hamiltonian decompositions of graphs, directed graphs and hypergraphs,
\emph{Ann. Discrete Math.} {\bf 3} (1978), 21--28.

\bibitem{Bermond} J.C. Bermond, A. Germa, M.C. Heydemann and D. Sotteau, Hypergraphes 
hamiltoniens, in \emph{Probl\`emes combinatoires et th\'eorie des graphes (Colloq. Internat. CNRS, Univ. Orsay, Orsay, 1976)}, 
vol. 260 of Colloq. Internat. CNRS, Paris (1973), 39--43. 

\bibitem{egt} B.~Bollob\'as,
\emph{Extremal Graph Theory}, p167,
Academic Press, 1978.

\bibitem{Bol83} B.~Bollob\'as, Almost all regular graphs are Hamiltonian,
\emph{European J. Combinatorics} {\bf 4} (1983), 97--106. 

\bibitem{Bol84} B.~Bollob\'as, 
The evolution of sparse graphs,
\emph{Graph theory and combinatorics},  Academic Press, London (1984),  35--57.

\bibitem{bollobasmatch} B.~Bollob\'as, 
Complete matchings in random subgraphs of the cube,
\emph{Random Structures \& Algorithms} {\bf  1} (1990), 95--104.



\bibitem{BF85} B.~Bollob\'as and A.~Frieze,
On matchings and Hamiltonian cycles in random graphs.  
\emph{Random graphs '83} (Poznan, 1983),  North-Holland Math. Stud., 118, North-Holland, Amsterdam (1985), 23--46.

\bibitem{bondy90} A.~Bondy,
Small cycle double covers of graphs. in: Cycles and rays (Montreal, PQ, 1987),  
\emph{NATO Adv. Sci. Inst. Ser. C Math. Phys. Sci.} {\bf 301}, Kluwer Acad. Publ., Dordrecht (1990), 21--40.

\bibitem{bondy} A.~Bondy,
Basic graph theory: paths and circuits,
in \emph{Handbook of Combinatorics}, Vol. 1, Elsevier, Amsterdam  (1995), 3--110.

\bibitem{BHSch} E.~Bu{\ss}, H.~H\'an and M. Schacht, Minimum vertex degree conditions for loose Hamilton cycles in 3-uniform hypergraphs,
\emph{J. Combin. Theory B}, to appear.

\bibitem{BHP} D.~Bryant, D.~Horsley and W.~Pettersson, 
Cycle decompositions V: Complete graphs into cycles of arbitrary lengths, 
\emph{Proc. London Math. Soc.}, to appear.

\bibitem{CBK} P.~Ch{\^a}u, L.~DeBiasio and H.A.~Kierstead, Pos\'a's conjecture for graphs of order at least $2\times 10^8$,
\emph{Random Structures \& Algorithms} {\bf 39} (2011), 507--525.

\bibitem{1factorization} A.G.~Chetwynd and A.J.W.~Hilton, 
Regular graphs of high degree are 1-factorizable, 
\emph{Proc. London Math. Soc.} {\bf 50} (1985), 193--206.

\bibitem{CH} A.G.~Chetwynd  and A.J.W. Hilton, 
$1$-factorizing regular graphs of high degree---an improved bound, 
\emph{Discrete Math.} {\bf 75} (1989), 103--112.

\bibitem{CHM} D.~Christofides, J.~Hladk\'y and A.~M\'ath\'e, Hamilton cycles in dense vertex-transitive graphs, \emph{J. Combin. Theory B}, to appear.

\bibitem{CKO} D.~Christofides, D.~K\"uhn and D.~Osthus, Edge-disjoint
Hamilton cycles in graphs, \emph{J. Combin. Theory B} \textbf{102} (2012), 1035--1060.


\bibitem{Chvatal} V.~Chv\'atal, Tough graphs and Hamiltonian circuits, \emph{Discrete Math.}~\textbf{5} (1973), 215--228.

\bibitem{paper2} B.~Csaba, D.~K\"uhn, A.~Lo, D.~Osthus and A.~Treglown, Proof of the $1$-factorization and Hamilton decomposition conjectures II: the bipartite case, preprint.

\bibitem{paper3} B.~Csaba, D.~K\"uhn, A.~Lo, D.~Osthus and A.~Treglown, Proof of the $1$-factorization and Hamilton decomposition conjectures III: approximate decompositions, preprint.



\bibitem{CucklerKahn} B.~Cuckler and J.~Kahn, Hamiltonian cycles in Dirac graphs, \emph{Combinatorica} {\bf 29} (2009), 299--326.

\bibitem{CucklerKahn2} B.~Cuckler and J.~Kahn, Entropy bounds for perfect matchings and Hamiltonian cycles, \emph{Combinatorica} {\bf 29} (2009), 327--335.

\bibitem{CM} A.~Czygrinow and T.~Molla, Tight co-degree condition for the existence of loose Hamilton cycles in 3-graphs, preprint.

\bibitem{Dirac} G.~Dirac, Some theorems on abstract graphs, \emph{Proc. London Math. Soc.} \textbf{2} (1952), 69--81.


\bibitem{DFloose} A.~Dudek and A. Frieze, Loose Hamilton cycles in random uniform hypergraphs,
\emph{Electronic J. Combinatorics} {\bf 18} (2011), \# P48.

\bibitem{DFtight} A.~Dudek and A. Frieze, Tight Hamilton cycles in random uniform hypergraphs,
\emph{Random Structures \& Algorithms} {\bf 42} (2013), 374--385.

\bibitem{DFLS} A.~Dudek, A. Frieze, P. Loh and S. Speiss, Optimal divisibility conditions for loose Hamilton cycles in random hypergraphs,
\emph{Electronic J. Combinatorics} {\bf 19} (2012), \# P44.


\bibitem{DFRS} A.~Dudek, A. Frieze, A. Ruci\'nski and M. {\v S}ileikis, Loose Hamilton cycles in regular hypergraphs,
\emph{Combin. Probab. Comput.}, to appear.

\bibitem{DFJ} M.~Dyer, A.~Frieze and M.~Jerrum,
Approximately counting Hamilton paths and cycles in dense graphs, \emph{SIAM J. Computing} {\bf 27} (1998), 1262--1272. 

\bibitem{FKS} A.~Ferber, M.~Krivelevich and B.~Sudakov,
Counting and packing Hamilton cycles in dense graphs and oriented graphs,  preprint.


\bibitem{Frcount} A.~Frieze,
On the number of perfect matchings and Hamilton cycles in $\eps$-regular non-bipartite graphs,
\emph{Electronic J. Combinatorics} {\bf 7} (2000), \# R57.

\bibitem{Frieze_loose} A.~Frieze, Loose Hamilton cycles in random 3-uniform hypergraphs,
\emph{Electronic J. Combinatorics} {\bf 17} (2010), \# N28.

\bibitem{fk} A.~Frieze and M.~Krivelevich, On packing Hamilton cycles in $\eps$-regular graphs,
\emph{J. Combin. Theory B}~{\bf{94}} (2005), 159--172.

\bibitem{FK08} A.~Frieze and M.~Krivelevich,
On two Hamilton cycle problems in random graphs,
\emph{Israel J. Math.}  {\bf 166} (2008), 221--234.   

\bibitem{FK} A. Frieze and M. Krivelevich, Packing Hamilton cycles in random and pseudo-random hypergraphs,
{\em Random Structures \& Algorithms} {\bf 41} (2012), 1--22.  

\bibitem{FKL}  A.~Frieze, M.~Krivelevich and P.~Loh,
Packing tight Hamilton cycles in 3-uniform hypergraphs,
\emph{Random Structures \& Algorithms} {\bf 40} (2012), 269--300.

\bibitem{GK} R.~Glebov and M.~Krivelevich, On the number of Hamilton cycles in sparse random graphs,
\emph{SIAM J. Discrete Math.} {\bf 27} (2013), 27--42. 

\bibitem{GKS} R.~Glebov, M.~Krivelevich and T.~Szab\'o, 
On covering expander graphs by Hamilton cycles, \emph{Random Structures \& Algorithms} (to appear).

\bibitem{gouldsurvey} R.~Gould, Advances on the Hamiltonian problem: a survey,
\emph{Graphs and Combinatorics} \textbf{19} (2003), 7--52.

\bibitem{gouldsurveyIII} R.~Gould, Recent advances on the Hamiltonian problem: Survey III,
\emph{Graphs and Combinatorics} \textbf{30} (2014), 1--46.

\bibitem{HS} H. H\`an and M. Schacht, Dirac-type results for loose Hamilton 
cycles in uniform hypergraphs, {\em J. Combin. Theory B} {\bf 100} (2010), 332--346.

\bibitem{HZ} J.~Han and Y.~Zhao, Minimum vertex degree threshold for loose Hamilton cycles in 3-uniform hypergraphs, preprint.

\bibitem{Hartkefactors} S.G.~Hartke, R.~Martin and T.~Seacrest,
Relating minimum degree and the existence of a $k$-factor, research
manuscript.

\bibitem{HartkeHCs} S.G.~Hartke and T.~Seacrest, Random partitions
and edge-disjoint Hamiltonian cycles, preprint.

\bibitem{HKS} D.~Hefetz, M.~Krivelevich and T.~Szab\'o, Hamilton cycles in highly connected and expanding graphs,
\emph{Combinatorica}~\textbf{29} (2009), 547--568. 

\bibitem{HKLO} D.~Hefetz, D.~K\"uhn, J.~Lapinskas and D.~Osthus,
Optimal covers with Hamilton cycles in random graphs, \emph{Combinatorica}, to appear.

\bibitem{HLS} H.~Huang, C.~Lee and B.~Sudakov, Bandwidth theorem for random graphs,
\emph{J. Combin. Theory B} {\bf 102} (2012), 14--37.

\bibitem{Jackson79} B. Jackson, Edge-disjoint Hamilton cycles in regular graphs of large degree,
\emph{J. London Math. Soc.} {\bf 19} (1979), 13--16.

\bibitem{jackson_for_Haggvist} B.~Jackson,
Hamilton cycles in regular 2-connected graphs,
\emph{J.~Combin. Theory B} {\bf 29} (1980), 27--46.

\bibitem{jackson} B.~Jackson, Long paths and cycles in oriented graphs, \emph{J. Graph Theory}~\textbf{5} (1981), 245--252.

\bibitem{jlz} B.~Jackson, H.~Li and Y.~Zhu,
Dominating cycles in regular 3-connected graphs,
\emph{Discrete Math.} {\bf 102} (1991), 163--176.

\bibitem{Janson1} S.~Janson, The numbers of spanning trees, Hamilton cycles and perfect matchings in a random
graph, \emph{Combin. Probab. Comput.} {\bf 3} (1994), 97--126.

\bibitem{Janson2} S.~Janson, Random regular graphs: asymptotic distributions and contiguity, \emph{Combin. Probab. Comput.} {\bf 4} (1995), 369--405.

\bibitem{JKV08} A.~Johansson, J.~Kahn and V.~Vu, Factors in random graphs, \emph{Random Structures \& Algorithms} {\bf 33} (2008), 1--28. 

\bibitem{jung} H.A.~Jung,
Longest circuits in 3-connected graphs,
\emph{Finite and infinite sets, Vol I, II, Colloq. Math. Soc. J\'anos Bolyai} {\bf 37} (1984), 403--438.

\bibitem{Karp} R.~Karp, Reducibility among combinatorial problems, Complexity of computer computations,
Plenum, New York, 1972, 85--103.

\bibitem{KK} G.Y.~Katona and H.A.~Kierstead, Hamiltonian chains in hypergraphs, \emph{J. Graph
Theory} {\bf 30} (1999), 205--212.

\bibitem{KKMO} P. Keevash, D. K\"uhn, R. Mycroft and D. Osthus, Loose Hamilton 
cycles in hypergraphs, \emph{Discrete Math.} {\bf 311} (2011), 544--559.

\bibitem{KW01} J.H. Kim and N. Wormald, 
Random matchings which induce Hamilton cycles and Hamiltonian decompositions of random regular graphs, 
\emph{J. Combin. Theory B} {\bf 81} (2001), 20--44. 

\bibitem{AHDoRG} F.~Knox, D.~K\"uhn and D.~Osthus,
Approximate Hamilton decompositions of random graphs, 
\emph{Random Structures \& Algorithms}~{\bf 40} (2012), 133--149.

\bibitem{KKO} F.~Knox, D.~K\"uhn and D.~Osthus, 
Edge-disjoint Hamilton cycles in random graphs, \emph{Random Structures \& Algorithms}, to appear.

\bibitem{KSSz98} J.~Koml\'os, G.~N.~S\'ark\"ozy and E.~Szemer\'edi,
Proof of the Seymour conjecture for large graphs,
\emph{Ann.~Combin.}~\textbf{2} (1998), 43--60.

\bibitem{KSz83} J.~Koml\'os and E.~Szemer\'edi, Limit distribution for the existence of Hamilton cycles in random
graphs, \emph{Discrete Math.} {\bf 43} (1983), 55--63.

\bibitem{Korshunov} A.D.~Korshunov, Solution of a problem of P. Erd\H{o}s and A. R\'enyi on Hamiltonian cycles in non-oriented graphs,
\emph{Diskret. Analiz.} {\bf 31} (1977), 17--56 (in Russian).


\bibitem{Kotzig}A. Kotzig,
Hamilton graphs and Hamilton circuits, \emph{Theory of Graphs and its Applications,
Proceedings of the Symposium of Smolenice, 1963}, Nakl. \v{C}SAV
Praha, (1964), 62--82.


\bibitem{Krivelevich} M.~Krivelevich, On the number of Hamilton cycles in pseudo-random graphs, \emph{Electronic J. Combinatorics} {\bf 19} (2012), \#P25.

\bibitem{KLeeS} M.~Krivelevich, C.~Lee and B.~Sudakov, Robust Hamiltonicity of Dirac graphs, \emph{Transactions Amer. Math. Soc.} {\bf 366} (2014), 3095--3130.


\bibitem{KrS} M. Krivelevich and W. Samotij, 
Optimal packings of Hamilton cycles in sparse random graphs, \emph{SIAM J.~Discrete Math.}~{\bf 26} (2012), 964--982.

\bibitem{KS} M.~Krivelevich and B.~Sudakov,  Sparse pseudo-random graphs are Hamiltonian,
\emph{J.~Graph Theory}~\textbf{42} (2003), 17--33.

\bibitem{KLOmindeg} D.~K\"uhn, J.~Lapinskas and D.~Osthus, Optimal packings of Hamilton cycles in graphs of high minimum degree,
\emph{Combin. Probab. Comput.} {\bf  22} (2013), 394--416.

\bibitem{KLOP} D.~K\"uhn, J.~Lapinskas, D.~Osthus and V.~Patel,
Proof of a conjecture of Thomassen on Hamilton cycles in highly connected tournaments, \emph{Proc. London Math. Soc}, to appear.

\bibitem{paper4} D.~K\"uhn, A.~Lo, and D.~Osthus, Proof of the $1$-factorization and Hamilton decomposition conjectures IV: exceptional systems for the two cliques case, preprint.

\bibitem{paper1} D.~K\"uhn, A.~Lo, D.~Osthus and A.~Treglown,
Proof of the $1$-factorization conjecture I: the $2$-clique case, preprint.

\bibitem{KLOS1} D. K\"uhn, A.~Lo, D.~Osthus and K.~Staden, The robust component structure of dense regular graphs and applications,
preprint.

\bibitem{KLOS2} D. K\"uhn, A.~Lo, D.~Osthus and K.~Staden, Solution to a problem of Bollob\'as and H\"aggkvist on Hamilton cycles in regular graphs,
preprint.

\bibitem{KMO} D. K\"uhn, R.~Mycroft and D.~Osthus,
Hamilton $\ell$-cycles in uniform hypergraphs, \emph{J. Combin. Theory A} {\bf 117} (2010), 910--927. 

\bibitem{KO} D.~K\"uhn and D.~Osthus, Loose Hamilton cycles in 3-uniform 
hypergraphs of high minimum degree, {\em J. Combin. Theory B} {\bf 96} 
(2006), 767--821. 

\bibitem{KOsurvey} D. K\"uhn and D. Osthus, A survey on Hamilton cycles in directed graphs, \emph{European J. Combinatorics}~\textbf{33}
(2012), 750--766.

\bibitem{KOrandposa} D. K\"uhn and D. Osthus, On Posa's conjecture for random graphs,
\emph{SIAM J. Discrete Math.} {\bf 26} (2012), 1440--1457.

\bibitem{Kelly} D.~K\"uhn and D.~Osthus, Hamilton decompositions of regular expanders: a proof of Kelly's conjecture for large tournaments,
\emph{Adv. in Math.} {\bf 237} (2013), 62--146.

\bibitem{KellyII} D.~K\"uhn and D.~Osthus,
Hamilton decompositions of regular expanders: applications, \emph{J. Combin. Theory B} {\bf 104} (2014), 1--27.

\bibitem{KOBerge} D.~K\"uhn and D.~Osthus, Decompositions of complete uniform hypergraphs into Hamilton Berge cycles, \emph{J. Combin. Theory A} {\bf 126} (2014), 128--135.

\bibitem{KOTchvatal} D.~K\"uhn, D.~Osthus and A.~Treglown,
Hamiltonian degree sequences in digraphs, 
\emph{J.~Combin. Theory  B} {\bf 100} (2010), 367--380.

\bibitem{KOTkelly} D.~K\"uhn, D.~Osthus and A.~Treglown,
Hamilton decompositions of regular tournaments, 
\emph{Proc.~London Math.~Soc.} {\bf 101} (2010), 303--335.


\bibitem{Kutnar} K.~Kutnar and D. Maru\v{s}i\v{c}, Hamilton cycles and paths in vertex-transitive graphs -- current directions, \emph{Discrete Math.}~\textbf{309}
(2009), 5491--5500.

\bibitem{LS} C.~Lee and B.~Sudakov, 
Dirac's theorem for random graphs, \emph{Random Structures \& Algorithms}~\textbf{41} (2012), 293--305.



\bibitem{Lovasz} L.~Lov\'asz. Problem 11. In 
\emph{Combinatorial Structures and
their Applications, Proceedings of the Calgary International Conference on Combinatorial Structures
and their Applications, R. Guy, H. Hanani, N. Sauer, and J. Sch\"onheim, editors},
Gordon and Breach Science Publishers, New York, 1970.

\bibitem{lucas} E.~Lucas, \emph{R\'ecr\'eations Math\'ematiques}, Vol.~2, Gautheir-Villars, 1892.

\bibitem{Marusic} D.~Maru\v{s}ic, Hamiltonian circuits in Cayley graphs, \emph{Discrete Math.} {\bf 46} (1983), 49--54.

\bibitem{moon} J.W.~Moon, \emph{Topics on tournaments}, 
Holt, Rinehart and Winston, New York, 1968.

\bibitem{MPW} T.~M\"uller, X.~P\'erez-Gimenez and N.~Wormald,
Disjoint Hamilton cycles in the random geometric graph,
\emph{J.~Graph Theory} {\bf 68} (2011), 299--322.

\bibitem{MP}  A.~Muthusamy and P.~Paulraja, Hamilton cycle decomposition of line graphs and a conjecture of Bermond,
{\em J. Combin. Theory B} {\bf 64} (1995), 1--16. 
   
\bibitem{NWreg} C.St.J.A.~Nash-Williams, 
Valency sequences which force graphs to have Hamiltonian circuits,
University of Waterloo Research Report, Waterloo, Ontario, 1969.

\bibitem{initconj} C.St.J.A.~Nash-Williams, Hamiltonian lines
in graphs whose vertices have sufficiently large valencies, in \emph{Combinatorial
theory and its applications, III (Proc. Colloq., Balatonf\"ured, 1969)},
North-Holland, Amsterdam (1970), 813--819.


\bibitem{Diracext} C.St.J.A.~Nash-Williams, Edge-disjoint
Hamiltonian circuits in graphs with vertices of large valency, in
\emph{Studies in Pure Mathematics (Presented to Richard Rado)}, Academic
Press, London (1971), 157--183.

\bibitem{OS} D.~Osthus and K.~Staden,
Approximate Hamilton decompositions of regular robustly expanding digraphs, 
\emph{SIAM J. Discrete Math.} {\bf 27} (2013), 1372--1409.


\bibitem{reed} L.~Perkovic and B.~Reed, 
Edge coloring regular graphs of high degree, \emph{Discrete Math.} {\bf 165/166} (1997), 567--578.

\bibitem{riordan} O.~Riordan,
Spanning subgraphs of random graphs,
\emph{Combin. Probab. Comput.} {\bf 9} (2000), 125--148.

\bibitem{RR} V. R\"odl and A.~Ruci\'nski, Dirac-type questions for hypergraphs -- a survey (or more problems for Endre to solve),
in: An Irregular Mind (Szemer\'edi is 70), \emph{Bolyai Soc. Math. Studies} {\bf 21} (2010), 561--590.

\bibitem{RR2} V. R\"odl and A.~Ruci\'nski, Families of triples with high minimum degree are Hamiltonian,
\emph{Discuss. Math. -- Graph Th.} {\bf 34} (2014) 363--383.

\bibitem{RRS} V. R\"odl, A.~Ruci\'nski and E. Szemer\'edi, A Dirac-type theorem 
for 3-uniform hypergraphs, {\em Combin. Probab. Comput.} {\bf 15} (2006), 
229--251.

\bibitem{RRS2} V. R\"odl, A. Ruci\'nski and E. Szemer\'edi, An approximate 
Dirac-type theorem for $k$-uniform hypergraphs, {\em Combinatorica} {\bf 28} 
(2008), 229--260.

\bibitem{RRS3} V. R\"odl, A. Ruci\'nski and E. Szemer\'edi, Dirac-type conditions for hamiltonian paths and cycles in 3-uniform hypergraphs,
\emph{Adv. in Math.} {\bf 227} (2011), 1225--1299. 

\bibitem{SSS} G.~S\'ark\"ozy, S.~Selkow and E.~Szemer\'edi, On the number of Hamiltonian cycles in Dirac
graphs, \emph{Discrete Math.} {\bf 265} (2003), 237--250.

\bibitem{schroeder} M.W.~Schroeder,
On Hamilton cycle decompositions of $r$-uniform $r$-partite hypergraphs,
{\em Discrete Math.} {\bf 315/316} (2014), 1--8.

\bibitem{stiebitz} M.~Stiebitz, D.~Scheide, B.~Toft and L.M.~Favrholdt,
\emph{Graph Edge Coloring: Vizing's Theorem and Goldberg's Conjecture}, Wiley 2012.

\bibitem{SV} B. Sudakov and V. Vu, Local resilience of graphs, \emph{Random Structures \& Algorithms} {\bf 33} (2008), 
409--433.


\bibitem{thomassenconj} C.~Thomassen, Edge-disjoint Hamiltonian paths and cycles in tournaments,
\emph{Proc. London Math. Soc.}~{\bf 45} (1982), 151--168.

\bibitem{vaughan} E.~Vaughan,
An asymptotic version of the multigraph $1$-factorization conjecture, 
\emph{J.~Graph Theory} {\bf 72} (2013),  19--29.

\bibitem{verrall} H.~Verrall, Hamilton decompositions of complete 3-uniform hypergraphs,
\emph{Discrete Math.} {\bf 132} (1994), 333--348.

\bibitem{wagner} D.~Wagner,
On the perfect $1$-factorization conjecture,
\emph{Discrete Math.} {\bf 104} (1992), 211--215.


\end{thebibliography}
\end{document}